\documentclass[12pt,reqno]{amsart} 
\usepackage{amsthm}

\usepackage{amsmath,amssymb}
\usepackage{ascmac}
\usepackage{mathrsfs}

\usepackage[abbrev]{amsrefs}

\usepackage{graphicx,xcolor}

\usepackage[%
 dvipdfmx,%
 setpagesize=false,%
 bookmarks=true,%
 bookmarksdepth=tocdepth,%
 bookmarksnumbered=true,%
 colorlinks=false,%
 pdftitle={},%
 pdfsubject={},%
 pdfauthor={},%
 pdfkeywords={}%
]{hyperref}            
\allowdisplaybreaks[1]
\theoremstyle{definition}
\newtheorem{dfn}{Definition}
\newtheorem{thm}{Theorem}[section]

\newtheorem{prop}[thm]{Proposition}
\newtheorem{lem}[thm]{Lemma}
\newtheorem{rem}{Remark}[section]

\newtheorem*{nota}{Notation}

\newcommand{\eqsp}[1]{{\begin{equation}\begin{split}#1\end{split}\end{
equation}}}

\makeatletter
\newcommand{\subscripts}[3]{%
  \@mathmeasure\z@\displaystyle{#2}%
  \global\setbox\@ne\vbox to\ht\z@{}\dp\@ne\dp\z@
  \setbox\tw@\box\@ne
  \@mathmeasure4\displaystyle{\copy\tw@_{#1}}%
  \@mathmeasure6\displaystyle{{#2}_{#3}}%
  \dimen@-\wd6 \advance\dimen@\wd4 \advance\dimen@\wd\z@
  \hbox to\dimen@{}\mathop{\kern-\dimen@\box4\box6}%
}
\makeatother

\begin{document}

\begin{center}\LARGE \bf
 On the uniqueness of the quasi-geostrophic equation with the fractional Laplacian
\end{center}
   
\footnote[0]
{


E-mail: *t-iwabuchi@tohoku.ac.jp, **okazaki.taiki.r5@dc.tohoku.ac.jp

}
\vskip5mm

\begin{center}
 {\large Tsukasa Iwabuchi* and Taiki Okazaki**} 
 \vskip2mm
 {\large Mathematical Institute, Tohoku University}\\
 {\large Sendai 980-8578 Japan}
\end{center}

\vskip5mm

\begin{center}
\begin{minipage}{135mm}
\footnotesize
{{\sc Abstract.} We consider the uniqueness of the solution of the surface quasi-geostrophic equation with fractional Laplacian. 
We show that the uniqueness holds in non-homogeneous Besov spaces without any additional assumption which is supposed to constract solutions. 
When the power of the fractional Laplacian is close to 2, 
we prove that the uniqueness with the regularity index $s=-1/2$. 
We extract the least regularity $s=-1/2$ for the well-definedness of the nonlinear term of the equation. }

\end{minipage}
\end{center}

\section{Introduction}
In this paper, we consider the quasi-geostrophic equation in $\mathbb{R}^2$.

\begin{equation}
   \label{SQG}
   \begin{cases}
     \partial_t \theta + \Lambda^\alpha \theta + (u \cdot \nabla) \theta = 0,
           &\quad  \   t>0,x\in{\mathbb{R}}^2, \\
     u = \nabla ^\perp \Lambda^{-1}\theta,
           &\quad  \ t>0,x\in{\mathbb{R}}^2,\\
     \theta(0,x) = \theta_0(x),
           &\quad  \ x\in{\mathbb{R}}^2,\\
   \end{cases}
\end{equation}
where $\Lambda^\alpha=(-\Delta)^{\alpha/2}$, $0<\alpha\leq 2$ and $\nabla ^\perp = (-\partial_{x_2},\partial_{x_1})$. 
The real-valued function $\theta(t,x)$ denotes the potential temperature 
and $u$ denotes the velocity of the fluid. 
The quasi-geostrophic equation is an important model in geophysical dynamics, 
which describes large scale atmospheric and oceantic motion with small Rossby numbers in the special case of constant
potential vorticity and buoyancy frequency (see \cite{Pe_1979}).
From a mathematical point of view, the quasi-geostrophic equation has similar stracture with the there-dimensional Euler and Navier-Stokes equations (see \cite{Co_Ma_Ta_1994}).

Let us recall several known result. 
The cases $\alpha>1$, $\alpha=1$, $\alpha<1$ are called sub-critical case, critical case, super-critical case, respectively. 
In sub-critical case, the global-in-time regularity for large initial data is obtained by Constantin and Wu \cite{Co_Wu_1999}. 
Also, in critical case, Constantin, Cordoba and Wu \cite{Co_Co_Wu_2001} proved the global existence and regularity for small initial data and 
Kiselev, Nazarov and Volberg \cite{Ki_Na_Vo_2007}, Caffarelli and Vasseur \cite{Ca_Va_2010} solved for large initial data case. 
On the other hand, in super-critical case, global-in-time regularity for large initial data is an open problem. 
The regularity for small initial data is known (see e.g. \cite{Co_Vi_2016}). 
In the framework of Besov spaces, Wang and Zhang \cite{Wa_Zh_2011} proved the existence of the local unique solution 
$\theta\in C([0,T];\mathring{B}^{1-\alpha}_{\infty,q})\cap \widetilde{L}^1(0,T; B^1_{\infty,q})$ for $\theta_0\in\mathring{B}^{1-\alpha}_{\infty,q},1\leq q<\infty$. 
Here $\mathring{B}^{1-\alpha}_{\infty,q}$ is the closure of the Schwartz class with respect to the norm of the non-homogenous Besov space $B^{1-\alpha}_{\infty,q}$.

Regarding the uniqueness of the solution, 
if $1< \alpha < 2$, then Ferreira \cite{Fe_2011} established that the uniqueness of the solution holds in $C([0,T];L^{2/(\alpha-1)})$. 
Also, for $\alpha=2$, Iwabuchi and Ueda \cite{I_U_2024} showed the uniqueness of solution $\theta\in C([0,T]; L^2)$. 
In critical case, the uniqueness of the solution in $C([0,T];\dot{B}_{\infty,1}^0)$ was proven in \cite{I_O_2025}. 
Recently, if $0<\alpha<3/2$, then Buckmaster, Shkoller and Vicol \cite{Bu_Sh_Vi_2019} proved that 
the nonuniqueness of the weak solution in $L_{\text{loc}}^2(0,T;\dot{H}^{-1/2}(\mathbb{T}^2))$. 

We also mention some literature on the Navier-Stokes equations. 
The scale critical spaces of the two-dimensional surface quasi-geostrophic equation are 
same as those of the two-dimensional incompressible Navier-Stokes equations. 
It is a classical result that the weak solution $u\in L^\infty(0,T;L^2)\cap L^2(0,T;\dot{H}^1)$ is unique if $u_0\in L^2$ for
the two-dimensional incompressible Navier-Stokes equations. 
The non-uniqueness in $C([0,T];L^p(\mathbb{T}^2))$ with $1\leq p<2$ for the two-dimensional incompressible Navier-Stokes equations
is provided by Cheskidov and Luo \cite{Ch_Lu_2023}. 
The uniqueness in the critical space $C([0,T];L^2)$ does not seem to be resolved. 
In the case when  the space dimension is three, 
Meyer \cite{Me_1996}, Furioli, Lemari\'{e}-Rieusset and Terraneo \cite{Fu_Le_Te_1997} and Monniaux \cite{Mo_1999} proved that 
the uniqueness of the solution in the critical space $C([0,T];L^3)$ for the incompressible Navier-Stokes equations. 
Also Lions and Masmoudi \cite{Li_Ma_2001} proved that the uniqueness in $L^3$ using a dual problem. 
Regarding the nonuniqueness result, Buckmaster and Vicol \cite{Bu_Vi_2019} showed that the nonuniqueness in $C([0,T];L^2(\mathbb{T}^3))$. 

In this paper, in sub-critical case, we prove that the uniqueness of the solution of \eqref{SQG} 
in larger spaces than the existing results \cite{Fe_2011, I_U_2024}. 
It is known that the following inclusion relation holds 
between the Besov spaces and the Legesgue spaces: 
\begin{equation*}
  L^{2/(\alpha-1)}\hookrightarrow B_{2/(\alpha-1),\infty}^0
  \hookrightarrow B_{p,\infty}^{1-\alpha+2/p}
  \hookrightarrow B_{\infty,\infty}^{1-\alpha} 
\end{equation*}
for $p> 2/(\alpha-1)$. 
In particular, we consider the uniqueness in Besov spaces which the regularity indices is negative. 
Also, if $0<\alpha\leq 1$, the existence of the local unique solution 
$\theta\in L^\infty(0,T;B^{1-\alpha}_{\infty,\infty})\cap \widetilde{L}^1(0,T; B^1_{\infty,\infty})$ for $\theta_0\in\mathring{B}^{1-\alpha}_{\infty,\infty}$ 
with smallness condition is proved by \cite{Wa_Zh_2011}. 
The space $B^{1-\alpha}_{\infty,\infty}$ is the largest critical non-homogeneous Besov space of \eqref{SQG}. 
When the interpolation index of the Besov space is infinite, 
it is not guaranteed that the linear solution of \eqref{SQG} is continuous with respect to time. 
The following conditions are known for the linear solution to be continuous. 

\noindent {\bf Proposition. }
Let $0<\alpha\leq 2$, $s\in\mathbb{R}$, $1\leq p\leq \infty$ and $\theta_0\in B^{1-\alpha}_{\infty,\infty}$. 
Then 
  \begin{equation*}
    \lim_{t\to 0}e^{-t\Lambda^\alpha}\theta_0=\theta_0 \text{ in } B^{s}_{p,\infty}
    \text{ if and only if }
    \lim_{j\to \infty}2^{sj}{\|\phi_j*\theta_0\|}_{L^p}=0. 
  \end{equation*}
\noindent 
Here $\phi_j$ ($j\in\mathbb{N}$) is a Littlewood-Paley dyadic decomposition function. 
Hence, we consider the uniqueness in the subspace of $B^{1-\alpha}_{\infty,\infty}$ 
in which the linear solution of \eqref{SQG} is continuous. 
In addition, it is know that 
\begin{equation*}
  \dot{B}_{\infty,1}^0\hookrightarrow L^\infty 
  \hookrightarrow B_{\infty,\infty}^0. 
\end{equation*}
Thus, in critical case, we prove that the uniqueness of the solution 
in larger space than the existing result \cite{I_O_2025}. 

Let us recall the definition of non-homogeneous Besov spaces. 
We refer to the book by Bahouri, Chemin and Danchin~\cite{Ba_Ch_Da_2011}. 
\begin{dfn}
  Let $\{\psi\}\cup\{\phi_j\}_{j\in\mathbb{N}}\subset\mathcal{S}(\mathbb{R}^d)$ be such that 
  \begin{equation*}
    \text{supp }\widehat{\psi}\subset\{\xi\in\mathbb{R}^d\ | |\xi| \leq 4/3\}, 
  \end{equation*}
  \begin{equation*}
    \text{supp } \widehat{\phi_j}\subset\{\xi\in\mathbb{R}^d\ |\ 3/4\cdot2^{j-1}\leq |\xi| \leq 8/3\cdot2^{j-1}\} 
    \text{ for any }j\in\mathbb{N}
  \end{equation*}
  and 
  \begin{equation*}
    \widehat{\psi}(\xi)+\sum_{j\in\mathbb{N}}\widehat{\phi_j}(\xi)=1 \text{ for any }\xi\in\mathbb{R}^d.   \end{equation*}
  For $s\in\mathbb{R}$\ and $1\leq p,q \leq \infty$,\ we define the non-homogeneous Besov spaces as follows. 
  \begin{equation*}
    B^s_{p,q}=B^s_{p,q}(\mathbb{R}^d):=\{f\in \mathcal{S}'(\mathbb{R}^d)\ |\ {\|f\|}_{B^s_{p,q}}<\infty\}, 
  \end{equation*}
  where
  \begin{equation*}
    {\|f\|}_{B^s_{p,q}}:={\|\psi*f\|}_{L^p}+{\left\|\left\{2^{sj}{\|\phi_j*f\|}_{L^p}\right\}_{j\in\mathbb{N}}\right\|}_{l^q}. 
  \end{equation*}
\end{dfn}

In this paper, we want to consider the product $(\nabla^\perp\Lambda^{-1}f)g$ for $f,g$ in the non-homogeneous Besov spaces 
with negative regularity indices. 
However, it is known that if $s<0$, then there exist $f,g\in B_{p,q}^s$ such that 
\begin{equation*}
  f g \notin \mathcal{S}'. 
\end{equation*}
On the other hand, the sum of two products
\begin{equation*}
  (\nabla^\perp\Lambda^{-1}f)g+(\nabla^\perp\Lambda^{-1}g)f
\end{equation*}
can be written in teams of second-order derivatives (see Lemma \ref{0325-5} and also \cite{I_U_2024}). 
In this case, if $-1/2<s< 0$, $2\leq p\leq \infty$ and $1\leq q\leq \infty$ 
or $s=-1/2$, $4\leq p\leq\infty$ and $1\leq q\leq p/(p-2)$, 
then for any $f,g\in B_{p,q}^s$, we can define $(\nabla^\perp\Lambda^{-1}f)g+(\nabla^\perp\Lambda^{-1}g)f$ in $\mathcal{S}'$. 
In the Appendix \ref{0403-2}, we show the definition of the sum of two product $(\nabla^\perp\Lambda^{-1}f)g+(\nabla^\perp\Lambda^{-1}g)f$ 
for $f,g$ in the non-homogeneous Besov spaces with the negative regularity indices 
and an example of the functions such that $(\nabla^\perp\Lambda^{-1}f)g\notin \mathcal{S}'$ when $s<0$. 
We refer to Lemma \ref{0423-1} for the justification of the sum of two product. 

The mild solution of \eqref{SQG} is defined as follows. 
\begin{dfn}\label{0328-1}
  Let $T>0$, $0<\alpha\leq 2$, $1\leq p,q\leq\infty$ and 
  $\theta_0\in B_{p,q}^{1-\alpha+2/p}$. 
  If $\theta:[0,T]\times\mathbb{R}^2\to\mathbb{R}$ satisfies 
  \begin{equation}\label{0327-3}
     \begin{cases}
      \theta \in C([0,T];B_{p,q}^{1-\alpha+2/p}), \\
      u\theta=(\nabla^\perp\Lambda^{-1}\theta)\theta \in\mathcal{S}' , \\
      \subscripts{\mathcal{S}'}{\displaystyle\big\langle\theta(t), \phi \big\rangle}{\mathcal{S}}
      =\subscripts{\mathcal{S}'}{\big\langle e^{-t\Lambda^\alpha}\theta_0,\phi\big\rangle}{\mathcal{S}}
      +\subscripts{\mathcal{S}'}{\left\langle \int_{0}^{t}e^{-(t-s)\Lambda^\alpha}(u\theta) ~{\rm d}s,\nabla\phi \right\rangle}{\mathcal{S}}\\ 
      \text{for all}\ t\in[0,T]\text{ and }\phi\in\mathcal{S}(\mathbb{R}^2),
     \end{cases}
  \end{equation}
  then we call $\theta$ a mild solution of \eqref{SQG}. 
\end{dfn}

We define the uniqueness of the mild solution. 
\begin{dfn}
  Let $T>0$ and $\theta$ is a mild solution of \eqref{SQG} for $B_{p,q}^{1-\alpha+2/p}$. 
  For any $\bar{\theta}$ which is a mild solution of \eqref{SQG} with $\bar{\theta}(0)=\theta(0)$, 
  if $\theta$ and $\bar{\theta}$ satisfies $\theta(t)=\bar{\theta}(t)$ in $B_{p,q}^{1-\alpha+2/p}$ for all $t\in [0,T]$, 
  then we say the mild solution of \eqref{SQG} in $C([0,T];B_{p,q}^{1-\alpha+2/p})$ is unique. 
\end{dfn}

Before stating our main result, we mention regularity index of $B_{p,q}^{1-\alpha+2/p}$. 
If $1-\alpha+2/p\geq -1/2$, then we can define $u\theta$. 
However, it is clear that 
\begin{equation*}
  1-\alpha+\frac{2}{p}\geq -\frac{1}{2} \text{ if and only if }
  \begin{cases}
    p\leq \frac{4}{2\alpha-3} & \frac{3}{2}< \alpha\leq 2\\ 
    p\leq \infty & 0<\alpha\leq \frac{3}{2}. 
  \end{cases}
\end{equation*}
Also, in the case of $1-\alpha+2/p= -1/2$, interpolation index is restricted. 
In this paper, we call the lowest regularity index of the scale critical non-homogeneous Besov space where the uniqueness holds the End-point, 
and we state our main theorem in the end-point case. 

\begin{thm}[End-point case]\label{0308-1}
  Let $T>0$, $0<\alpha\leq 2$. 
  \begin{enumerate}
    \item \label{0308-2}If $3/2<\alpha\leq 2$ and $p=4/(2\alpha-3)$, 
    then the mild solution of \eqref{SQG} in $ C([0,T]; B_{p,p/(p-2)}^{-1/2})$ is unique. 
    \item If $\alpha=3/2$ and $\theta_0$ satisfies ${\|\nabla^\perp\Lambda^{-1}(\psi*\theta_0)\|}_{L^\infty}<\infty$, 
    then the mild solution of \eqref{SQG} in $C([0,T]; B_{\infty, 1}^{-1/2})$ is unique. 
    \item \label{0308-4}If $0<\alpha<3/2$ and $\theta_0$ satisfies 
    ${\|\nabla^\perp\Lambda^{-1}(\psi*\theta_0)\|}_{L^\infty}<\infty$ 
    and $\displaystyle\lim_{j\to\infty}2^{(1-\alpha)j}{\|\phi_j*\theta_0\|}_{L^\infty}=0$, 
    then the mild solution of \eqref{SQG} in $C([0,T]; B_{\infty, \infty}^{1-\alpha})$ is unique. 
  \end{enumerate}
\end{thm}

\begin{rem}
  For $\theta\in B_{\infty,q}^s$, the Riesz transform of $\theta_0$ does not belong to $B_{\infty,q}^s$ in general. 
  Thus, in our theorem, we impose ${\|\nabla^\perp\Lambda^{-1}(\psi*\theta_0)\|}_{L^\infty}<\infty$. 
\end{rem}

In the Theorem \ref{0308-1}, if $3/2\leq\alpha<2$, 
then there is no inclusion relationship between the space where we prove the uniqueness 
and $L^{2/(\alpha-1)}$. 
The uniqueness in $L^{2/(\alpha-1)}$ is proven in the existing results \cite{Fe_2011, I_U_2024}. 
However, in the non-end-point i.e. $1-\alpha+2/p>-1/2$, 
we can prove the following theorem. 

\begin{thm}[Non-end-point case]\label{0402-1}
  Let $T>0$, $3/2\leq \alpha\leq 2$. 
  \begin{enumerate}
    \item If $3/2<\alpha\leq 2$, $2/(\alpha-1)\leq p<4/(2\alpha-3)$ 
    and $\theta_0$ satisfies $\displaystyle\lim_{j\to\infty}2^{(1-\alpha+2/p)j}{\|\phi_j*\theta_0\|}_{L^p}=0$, 
    then the mild solution of \eqref{SQG} in $C([0,T];B_{p,\infty}^{1-\alpha+2/p})$ is unique. 
    \item If $\alpha=3/2$, $4\leq p<\infty$ 
    and $\theta_0$ satisfies $\displaystyle\lim_{j\to\infty}2^{(-1/2+2/p)j}{\|\phi_j*\theta_0\|}_{L^p}=0$, 
    then the mild solution of \eqref{SQG} in $C([0,T];B_{p,\infty}^{-1/2+2/p})$ is unique. 
  \end{enumerate}
\end{thm}

In this paper, we only show the proof of Theorem \ref{0308-1}. 
For Theorem \ref{0402-1}, it can be proven in a similar proof to Theorem \ref{0308-1} \eqref{0308-2}. 
In the proof of Theorem \ref{0308-1}, 
we use structure of the nonlinear term of \eqref{SQG} written by the second order derivatives, 
and if $1\leq \alpha\leq 3/2$, then we devide the solution of \eqref{SQG} into high frequency part and low frequency part. 
Moreover, if $0<\alpha\leq 1$, then we consider the norm $B_{\infty,\infty}^{-s}$ with some small positive number $s$ 
instead of the $B_{\infty,\infty}^{1-\alpha}$ norm. 
To this end, we introduce some notations. 
\begin{nota}
  For any $j\in\mathbb{N}$, we denote 
  \begin{equation}\label{0308-3}
    f_j:=\phi_j*f,\quad S_jf:=\psi*f+\sum_{k\leq j}\phi_k*f \text{ and }
    \widetilde{S}_jf:=f-S_jf. 
  \end{equation}
\end{nota}

\section{Preliminarys}
In this section, we introduce several lemmas which are elemental properties of the non-homogeneous Besov spaces and heat kernel. 

\begin{lem}[\cite{Ba_Ch_Da_2011}]\label{0325-2}
  Let $1\leq p\leq \infty$.
  Then there exists constant $C>0$ such that 
  for any $j\in\mathbb{N}$ and $f$ with $\phi_j*f\in L^p$, we have 
  \begin{equation*}
     {\|\nabla(\phi_j* f)\|}_{L^p}\leq C2^j{\|\phi_j*f\|}_{L^p},\ {\|\Lambda^{-1}(\phi_j* f)\|}_{L^p}\leq C2^{-j}{\|\phi_j*f\|}_{L^p}. 
  \end{equation*}
  Also, there exists constant $C>0$ such that 
  for any $f$ with $\psi*f\in L^p$, we have 
  \begin{equation*}
     {\|\nabla(\psi* f)\|}_{L^p}\leq C{\|\psi *f\|}_{L^p}. 
  \end{equation*}
\end{lem}

\begin{lem}[\cite{Wu_Yu_2008}]\label{0327-1}
  Let $t>0$, $0<\alpha\leq 2$, $1\leq p\leq \infty$ and $j\in\mathbb{N}$. 
  Then there exist constants $C,c>0$ such that for any $f$ with $\phi_j*f\in L^p$, we have 
  \begin{equation*}
    {\|e^{-t\Lambda^\alpha}(\phi_j*f)\|}_{L^p}\leq Ce^{-ct2^{\alpha j}}{\|\phi_j*f\|}_{L^p}. 
  \end{equation*}
\end{lem}

\begin{lem}[\cite{Ba_Ch_Da_2011}]\label{0325-4}
  Let $s\in\mathbb{R}$, $1\leq p_1\leq p_2\leq\infty$ and $1\leq q_1\leq q_2\leq\infty$. 
  Then $B_{p_1,q_1}^s(\mathbb{R}^d)\hookrightarrow B_{p_2,q_2}^{s-d(1/p_1-1/p_2)}(\mathbb{R}^d)$ 
  i.e. there exist constant $C>0$ such that 
  for any $f\in B_{p_1,q_1}^s(\mathbb{R}^d)$, we have 
  \begin{equation*}
    {\|f\|}_{B_{p_2,q_2}^{s-d(1/p_1-1/p_2)}}\leq C{\|f\|}_{B_{p_1,q_1}^s}. 
  \end{equation*}
\end{lem}

\begin{lem}[\cite{Ba_Ch_Da_2011}]\label{0404-5}
  Let $2\leq p\leq\infty$. 
  Then $L^p\hookrightarrow B_{p,p}^0$
  i.e. there exist constant $C>0$ such that 
  for any $f\in L^p$, we have 
  \begin{equation*}
    {\|f\|}_{B_{p,p}^0}\leq C{\|f\|}_{L^p}.
  \end{equation*}
\end{lem}

\begin{lem}[\cite{Ba_Ch_Da_2011}]\label{0326-5}
  Let $s,s'\in\mathbb{R}$, $1\leq p,q\leq\infty$ and 
  $\alpha$ be a multi-index. 
  If a Fourier multiplier $m\in C^\infty(\mathbb{R}^d)$ satisfies that 
  there exists constant $C_\alpha>0$ such that for any $\xi\in\mathbb{R}^d$, 
  \begin{equation*}
    |\partial^\alpha m(\xi)|\leq C_\alpha(1+|\xi|)^{s'-|\alpha|}, 
  \end{equation*}
  then there exists constant $C>0$ such that for any $f\in B_{p,q}^{s+s'}$, we have 
  \begin{equation*}
    {\left\|\mathcal{F}^{-1}[m(\xi)\widehat{f}(\xi)]\right\|}_{B_{p,q}^s}
    \leq C{\|f\|}_{B_{p,q}^{s+s'}}. 
  \end{equation*}
\end{lem}

\begin{rem}
  Let $s\in\mathbb{R}$ and $1\leq q\leq\infty$. 
  Due to Lemma \ref{0326-5}, there exists constant $C>0$ such that 
  for any $f\in B_{\infty,q}^s$, we have 
  \begin{equation*}
    {\|\nabla\nabla^\perp\Lambda^{-1}f\|}_{B_{\infty,q}^{s-1}}
    \leq C{\|f\|}_{B_{\infty,q}^s}. 
  \end{equation*}
\end{rem}

\begin{lem}[\cite{Ba_Ch_Da_2011}]\label{0325-3}
  Let $s\in\mathbb{R}$, $\epsilon>0$ and $1\leq p,q,q_1,q_2\leq\infty$ with $1/q= 1/q_1+1/q_2$.
  Then there exists constant $C>0$ such that 
  for any $f\in B^{-\epsilon}_{\infty,q_1}$ and $g\in B^s_{p,q}$, we have 
  \begin{equation*}
    {\left\|\sum_{l\geq 2}(\psi*f)(\phi_l*g)+\sum_{k\leq l-2}(\phi_k*f)(\phi_l*g)\right\|}_{B^{s-\epsilon}_{p,q}}\leq C{\|f\|}_{B^{-\epsilon}_{\infty,q_1}}{\|g\|}_{B^s_{p,q_2}}. 
  \end{equation*}
\end{lem}

The following lemma is the bilinear estimate in the non-homogeneous Besov spaces with negative regularity indices, 
which is used for our proof of Theorem \ref{0308-1}. 

\begin{lem}\label{0325-5}
  Let $s>-2$, $s'\in\mathbb{R}$, $2\leq p,p_1,p_2\leq\infty$ with $1/p=1/p_1+1/p_2$ 
  and $1\leq q,q_1,q_2\leq\infty$ with $1/q=1/q_1+1/q_2$. 
  Then there exists constant $C>0$ such that for any $f\in B_{p_1,q_1}^{s'}$ and $g\in B_{p_2,q_2}^{s+1-s'}$, 
  we have 
  \begin{align}\label{0326-1}
    &{\left\|\sum_{|k-l|\leq 1}\left(\left(\nabla^\perp\Lambda^{-1}(\phi_k*f)\right)\cdot\nabla(\phi_l*g)
    +\left(\nabla^\perp\Lambda^{-1}(\phi_l*g)\right)\cdot\nabla(\phi_k*f)\right)\right\|}_{B_{p,q}^{s}}\notag\\ 
    \leq& C{\|f\|}_{B_{p_1,q_1}^{s'}}{\|g\|}_{B_{p_2,q_2}^{s+1-s'}}. 
  \end{align}
  Also, if $1=1/q_1+1/q_2$, then there exists constant $C>0$ such that for any $f\in B_{p_1,q_1}^{s'}$ and $g\in B_{p_2,q_2}^{-1-s'}$, 
  we have 
  \begin{align}\label{0326-2}
    &{\left\|\sum_{|k-l|\leq 1}\left(\left(\nabla^\perp\Lambda^{-1}(\phi_k*f)\right)\cdot\nabla(\phi_l*g)
    +\left(\nabla^\perp\Lambda^{-1}(\phi_l*g)\right)\cdot\nabla(\phi_k*f)\right)\right\|}_{B_{p,p}^{-2}}\notag\\ 
    \leq& C{\|f\|}_{B_{p_1,q_1}^{s'}}{\|g\|}_{B_{p_2,q_2}^{-1-s'}}. 
  \end{align}
\end{lem}

\begin{rem}
  In general, the non-negativity $s\geq 0$ is essentially needed for the product estimate (see section 4.8 in \cite{Ru_Si_1996}). 
  On the other hund, in Lemma \ref{0325-5}, we can estimate the sum of two products $(\nabla^\perp\Lambda^{-1}f)\cdot\nabla g+(\nabla^\perp\Lambda^{-1}g)\cdot\nabla f$ 
  on $B_{p,q}^s$ with $s>-2$ or $B_{p,p}^{-2}$ 
  since $(\nabla^\perp\Lambda^{-1}f)\cdot\nabla g+(\nabla^\perp\Lambda^{-1}g)\cdot\nabla f$ can be written by the second order derivatives as in the proof below. 
\end{rem}

\begin{proof}[Proof of Lemma \ref{0325-5}]
  Since $\nabla\cdot\nabla^\perp f=0$, we have 
  \begin{align}\label{0404-4}
    &\big(\nabla^\perp\Lambda^{-1}(\phi_k*f)\big)\cdot\nabla(\phi_l*g)
    +\big(\nabla^\perp\Lambda^{-1}(\phi_l*g)\big)\cdot\nabla(\phi_k*f)\notag\\
    =&\nabla\cdot\Big(\big(\nabla^\perp\Lambda^{-1}(\phi_k*f)\big)(\phi_l*g)
    -\big(\Lambda^{-1}(\phi_l*g)\big)\big(\nabla^\perp(\phi_k*f)\big)\Big). 
  \end{align}
  First, we prove in the case $p<\infty$. 
  We can write \eqref{0404-4} using the Fourier multiplier (see \cite{I_U_2024}), 
  \begin{align}
    &\nabla\cdot\Big(\big(\nabla^\perp\Lambda^{-1}(\phi_k*f)\big)(\phi_l*g)
    -\big(\Lambda^{-1}(\phi_l*g)\big)\big(\nabla^\perp(\phi_k*f)\big)\Big)\notag\\ 
    =&\nabla\cdot\left(\nabla\cdot m(D_1,D_2)\left(\nabla^\perp\Lambda^{-1}(\phi_k*f),\Lambda^{-1}(\phi_l*g)\right)\right), 
  \end{align}
  where the Fourier multiplier $m(\cdot,\cdot)$ is bounded in $L^p$ (see \cite{Gr_Mi_To_2013}). 
  Thus, by Lemma \ref{0325-2}, we have 
  \begin{align*}
    &{\left\|\sum_{|k-l|\leq 1}\left(\left(\nabla^\perp\Lambda^{-1}(\phi_k*f)\right)\cdot\nabla(\phi_l*g)
    +\left(\nabla^\perp\Lambda^{-1}(\phi_l*g)\right)\cdot\nabla(\phi_k*f)\right)\right\|}_{B_{p,q}^{s}}\notag\\ 
    \leq &C\sum_{|k-l|\leq 1}{\left\|\nabla^\perp\Lambda^{-1}(\phi_k*f)\right\|}_{L^{p_1}}{\left\|\Lambda^{-1}(\phi_l*g)\right\|}_{L^{p_2}}\\ 
    &+C\left(\sum_{j\in\mathbb{N}}\left(\sum_{k\geq j-4}\sum_{|k-l|\leq 1}2^{(s+2)j}{\left\|\nabla^\perp\Lambda^{-1}(\phi_k*f)\right\|}_{L^{p_1}}{\left\|\Lambda^{-1}(\phi_l*g)\right\|}_{L^{p_2}}\right)^q\right)^{1/q}. 
  \end{align*}
  Since $s+2>0$, we obtain \eqref{0326-1}. 
  In the case of $s=-2$ and $q=1$, using Lemma \ref{0404-5}, we get 
  \begin{align*}
    &{\left\|\sum_{|k-l|\leq 1}\left(\left(\nabla^\perp\Lambda^{-1}(\phi_k*f)\right)\cdot\nabla(\phi_l*g)
    +\left(\nabla^\perp\Lambda^{-1}(\phi_l*g)\right)\cdot\nabla(\phi_k*f)\right)\right\|}_{B_{p,p}^{-2}}\notag\\ 
    \leq&C{\left\|\sum_{|k-l|\leq 1}m(D_1,D_2)\left(\nabla^\perp\Lambda^{-1}(\phi_k*f),\Lambda^{-1}(\phi_l*g)\right)\right\|}_{L^p}\notag\\ 
    \leq&C\sum_{|k-l|\leq 1}{\left\|\nabla^\perp\Lambda^{-1}(\phi_k*f)\right\|}_{L^{p_1}}{\left\|\Lambda^{-1}(\phi_l*g)\right\|}_{L^{p_2}}. 
  \end{align*}
  By H\"older's inequality, we obtain \eqref{0326-2}. 
  If $p=\infty$, then we consider the variables of $f$ and $g$ separately (see \cite{I_O_2025}). 
  We have 
  \begin{align*}
    &\big(\nabla^\perp\Lambda^{-1}(\phi_k*f)\big)(x)(\phi_l*g)(x')
    -\big(\Lambda^{-1}(\phi_l*g)\big)(x')\big(\nabla^\perp(\phi_k*f)\big)(x') \notag\\
    =&\nabla_{\mathbb{R}^4}\cdot \mathcal{F}^{-1}_{\mathbb{R}^4}\left[\frac{-i(-\xi,\eta)}{|\xi|+|\eta|}\frac{i\xi^\perp}{|\xi|}\widehat{\phi}_k(\xi)\widehat{f}(\xi)\frac{1}{|\eta|}\widehat{\phi}_l(\eta)\widehat{g}(\eta)\right](x,x').  
  \end{align*}
  where $\mathcal{F}^{-1}_{\mathbb{R}^4}$ is the Fourier inverse transform in $\mathbb{R}^4$, 
  and $\nabla_{\mathbb{R}^4}$ is the gradient in $\mathbb{R}^4$. 
  Hence, by a similar estimate to $p<\infty$, 
  we can prove Lemma \ref{0325-5}. 
\end{proof}

We introduce the Chemin-Lerner space $\widetilde{L}^r(0,T;B_{p,q}^s)$. 

\begin{dfn}
  Let $T>0$, $s\in\mathbb{R}$ and $1\leq p,q,r\leq\infty$. 
  We define the the Chemin-Lerner space as follows. 
  \begin{equation*}
    \widetilde{L}^r(0,T;B_{p,q}^s)=\widetilde{L}^r(0,T;B_{p,q}^s(\mathbb{R}^d))
    :=\{f\in \mathcal{D}'((0,T)\times\mathbb{R}^d)|\ {\|f\|}_{\widetilde{L}^r(0,T;B_{p,q}^s)}<\infty\}, 
  \end{equation*}
  where 
  \begin{equation*}
    {\|f\|}_{\widetilde{L}^r(0,T;B_{p,q}^s)}:={\|\psi*f\|}_{L^r(0,T;L^p)}
    +{\left\|\left\{2^{sj}{\|\phi_j*f\|}_{L^r(0,T;L^p)}\right\}_{j\in\mathbb{N}}\right\|}_{l^q(\mathbb{N})}. 
  \end{equation*}
\end{dfn}

Next Lemma follows from Minkowski's inequality. 
\begin{lem}
  Let $T>0$, $s\in\mathbb{R}$ and $1\leq p,q,r\leq\infty$. 
  If $r\leq q$, then there exists some constant $C>0$ such that 
  for any $f\in L^r(0,T;B_{p,q}^s)$, we have 
  \begin{equation*}
    {\|f\|}_{\widetilde{L}^r(0,T;B_{p,q}^s)}
    \leq C{\|f\|}_{L^r(0,T;B_{p,q}^s)}.
  \end{equation*}
  Also, if $r\geq q$, then there exists some constant $C>0$ such that 
  for any $f\in \widetilde{L}^r(0,T;B_{p,q}^s)$, we have 
  \begin{equation*}
    {\|f\|}_{L^r(0,T;B_{p,q}^s)}
    \leq C{\|f\|}_{\widetilde{L}^r(0,T;B_{p,q}^s)}.
  \end{equation*}
\end{lem}

We introduce the product estimate and the commutator estimate 
in the non-homogenous Besov spaces. 

\begin{lem}[\cite{Ch_Mi_Zh_2007}]\label{0325-8}
  Let $T>0$, $2\leq p\leq\infty$, $1\leq q,q_1,q_2,r,r_1,r_2\leq\infty$ with $1/q=1/q_1+1/q_2$ and $1/r=1/r_1+1/r_2$ 
  and $s,\epsilon\in \mathbb{R}$ with $s+2/p+1>0$, $\epsilon>0$ and $s+\epsilon<2/p$. 
  Then there exists constant $C>0$ such that for all $u\in \widetilde{L}^{r_1}(0,T;(B_{p,q_1}^{2/p-\epsilon}(\mathbb{R}^2))^2)$ 
  with $\nabla\cdot u=0$ and $\theta\in\widetilde{L}^{r_2}(0,T;B_{p,q_2}^{s+\epsilon+1}(\mathbb{R}^2))$, we have 
  \begin{equation*}
    {\|u\cdot\nabla\theta\|}_{\widetilde{L}^{r}(0,T;B_{p,q}^{s})}
    \leq C{\|u\|}_{\widetilde{L}^{r_1}(0,T;B_{p,q_1}^{2/p-\epsilon})}{\|\theta\|}_{\widetilde{L}^{r_2}(0,T;B_{p,q_2}^{s+\epsilon+1})}.
  \end{equation*}
\end{lem}

\begin{lem}[\cite{Ch_Mi_Zh_2007,Da_2003}]\label{0325-7}
  Let $T>0$, $2\leq p\leq \infty$, $1\leq q,q_1,q_2,r,r_1,r_2\leq\infty$ with $1/q=1/q_1+1/q_2$ and $1/r=1/r_1+1/r_2$ 
  and $s,\epsilon\in \mathbb{R}$ with $s+2/p+1>0$, $0<\epsilon<2/p+1$ and $s+\epsilon<2/p+1$. 
  Then there exists constant $C>0$ such that for all $u\in \widetilde{L}^{r_1}(0,T;(B_{p,q_1}^{2/p-\epsilon+1}(\mathbb{R}^2))^2)$ 
  with $\nabla\cdot u=0$ and $\theta\in\widetilde{L}^{r_2}(0,T;B_{p,q_2}^{s+\epsilon}(\mathbb{R}^2))$, we have 
  \begin{align*}
    &{\|[\psi*,u]\cdot\nabla\theta\|}_{L^r(0,T;L^p)}+{\|\{2^{sj}{\|[\phi_j*,u]\cdot\nabla\theta\|}_{L^r(0,T;L^p)}\}_{j\in\mathbb{N}}\|}_{l^q}\notag\\
    \leq& C{\|\nabla u\|}_{\widetilde{L}^{r_1}(0,T;B_{p,q_1}^{2/p-\epsilon})}{\|\nabla\theta\|}_{\widetilde{L}^{r_2}(0,T;B_{p,q_2}^{s+\epsilon-1})}.
  \end{align*}
\end{lem}

\begin{lem}[\cite{Wa_Zh_2011}]\label{0326-9}
  Let $j\in\mathbb{N}$. 
  Then there exists constant $C>0$ such that for any $f$ and $g$ with 
  $\psi*f,\phi_j*f,\psi*g, \phi_j*g\in L^\infty$, we have 
  \begin{align*}
    &{\|[\nabla^\perp\Lambda^{-1}\psi*, \nabla^\perp\Lambda^{-1}f]\cdot \nabla g\|}_{L^\infty}\\
    \leq& C \Bigg({\|\psi*g\|}_{L^\infty}\sum_{k>4}{\|\phi_k*f\|}_{L^\infty}
    +{\|\psi*f\|}_{L^\infty}{\|\psi*g\|}_{L^\infty}
    +{\|\psi*f\|}_{L^\infty}{\|\phi_1*g\|}_{L^\infty}\\ 
    &\quad+{\|\phi_1*f|}_{L^\infty}{\|\psi*g\|}_{L^\infty} 
    +\sum_{|k-l|\leq 6}{\|\phi_k*f\|}_{L^\infty}{\|\phi_l*g\|}_{L^\infty}\Bigg).
  \end{align*}
\end{lem}

By the smoothing effect of heat kernel (Lemma \ref{0327-1}) 
and Lemma \ref{0325-5}, the nonlinear team in a mild solution of \eqref{SQG} 
belongs to $L^p$. 

\begin{lem}\label{0327-2}
  Let $T>0$ and $0<\alpha\leq 2$. 
  \begin{enumerate}
    \item If $3/2<\alpha\leq 2$ and $p=4/(2\alpha-3)$, then there exists constant $C>0$ such that 
    for any $\theta\in C([0,T];B_{p,p/(p-2)}^{-1/2})$, we have 
    \begin{equation*}
      \sup_{t\in [0,T]}{\left\|\int_{0}^{t}e^{-(t-s)\Lambda^\alpha}\left((\nabla^\perp\Lambda^{-1}\theta)\theta\right)~{\rm d}s\right\|}_{L^p}
    \leq C\max\{T, T^{1/2\alpha}\}{\|\theta\|}_{L^\infty(0,T;B_{p,p/(p-2)}^{-1/2})}^2. 
    \end{equation*}
    \item If $\alpha=3/2$, then there exists constant $C>0$ such that 
    for any $\theta\in C([0,T];B_{\infty,1}^{-1/2})$ with ${\|\nabla^\perp\Lambda^{-1}(\psi*\theta)\|}_{L^\infty(0,T; L^\infty)}<\infty$, we have 
    \begin{align*}
      &\sup_{t\in [0,T]}{\left\|\int_{0}^{t}e^{-(t-s)\Lambda^{3/2}}\left((\nabla^\perp\Lambda^{-1}\theta)\theta\right)~{\rm d}s\right\|}_{L^\infty}\\ 
      \leq& C\max\{T, T^{1/3}\}{\|\nabla^\perp\Lambda^{-1}\theta\|}_{L^\infty(0,T;B_{\infty,1}^{-1/2})}{\|\theta\|}_{L^\infty(0,T;B_{\infty,1}^{-1/2})}.
    \end{align*}
    \item If $1<\alpha<3/2$, then there exists constant $C>0$ such that 
    for any $\theta\in C([0,T];B_{\infty,\infty}^{1-\alpha})$ with ${\|\nabla^\perp\Lambda^{-1}(\psi*\theta)\|}_{L^\infty(0,T; L^\infty)}<\infty$, we have 
    \begin{align*}
      &\sup_{t\in [0,T]}{\left\|\int_{0}^{t}e^{-(t-s)\Lambda^\alpha}\left((\nabla^\perp\Lambda^{-1}\theta)\theta\right)~{\rm d}s\right\|}_{L^\infty}\\ 
      \leq& C\max\{T, T^{2-2/\alpha}\}{\|\nabla^\perp\Lambda^{-1}\theta\|}_{L^\infty(0,T;B_{\infty,\infty}^{1-\alpha})}{\|\theta\|}_{L^\infty(0,T;B_{\infty,\infty}^{1-\alpha})}.
    \end{align*}
    \item If $0<\alpha\leq 1$, then there exists constant $C>0$ such that 
    for any $\theta\in C([0,T];B_{\infty,\infty}^{1-\alpha})$ with ${\|\nabla^\perp\Lambda^{-1}(\psi*\theta)\|}_{L^\infty(0,T; L^\infty)}<\infty$, we have 
    \begin{align*}
      &\sup_{t\in [0,T]}{\left\|\int_{0}^{t}e^{-(t-s)\Lambda^\alpha}\left((\nabla^\perp\Lambda^{-1}\theta)\theta\right)~{\rm d}s\right\|}_{L^\infty}\\ 
      \leq& C\max\{T, T^{1/2}\}{\|\nabla^\perp\Lambda^{-1}\theta\|}_{L^\infty(0,T;B_{\infty,\infty}^{1-\alpha})}{\|\theta\|}_{L^\infty(0,T;B_{\infty,\infty}^{1-\alpha})}.
    \end{align*}
  \end{enumerate}
\end{lem}

We can justify the differentiation of $\theta_j=\phi_j*\theta$ with respect to time. 

\begin{lem}\label{0301-33}
   Let $T>0$, $3/2<\alpha\leq 2$, $p=4/(2\alpha-3)$ and 
   $\theta\in C([0,T];B_{p,p/(p-2)}^{-1/2})$\ be a mild solution of \eqref{SQG}.
   Then for any $j\in\mathbb{N}$, $\partial_t\theta_j\in C([0,T];B_{p,p/(p-2)}^{-1/2})$ and for all $t\in(0,T)$, $\theta_j$ satisfies 
   \begin{equation*}
    \partial_t \theta_j + \Lambda \theta_j + \nabla \phi_j * (u\theta) = 0\ \text{in}\ B_{p,p/(p-2)}^{-1/2}.
   \end{equation*}
   Moreover, a similar statement holds for $\theta\in C([0,T];B_{\infty,1}^{-1/2})$ 
   or $\theta\in C([0,T];B_{\infty,\infty}^{1-\alpha})$ 
   that satisfies ${\|\nabla^\perp\Lambda^{-1}(\psi*\theta)\|}_{L^\infty(0,T; L^\infty)}<\infty$. 
\end{lem}
\begin{proof}
  By smoothing effect (Lemma \ref{0327-1}) and Lemma \ref{0327-2}, 
  the linear and nonlinear terms of a mild solution of \eqref{SQG} belong to $L^p$. 
  Therefore, we can write \eqref{0327-3} in an integral form, 
  \begin{equation*}
    \int_{\mathbb{R}^2}\theta(t,x)\phi(x) ~{\rm d}x
    =\int_{\mathbb{R}^2}e^{-t\Lambda^\alpha}\theta_0(x)\,\phi(x)~{\rm d}x
    +\int_{\mathbb{R}^2}\left(\int_{0}^{t}e^{-(t-s)\Lambda}(u\theta) ~{\rm d}s\right)\cdot\nabla\phi(x) ~{\rm d}x. 
  \end{equation*}
  Thus, by a similar argument to Proposition 2.1 in \cite{I_O_2025}, 
  we can prove Lemma \ref{0301-33}. 
\end{proof}

\begin{lem}[\cite{Co_Co_2004}]\label{0325-1}
  Let $0\leq \alpha\leq 2$ and $1\leq p\leq\infty$. 
  Let $\theta,u$ and $f$ be smooth functions on $(0,\infty)\times\mathbb{R}^2$ 
  such that $\nabla\cdot u=0$. 
  If $\theta$ satisfies $\partial_t\theta+\Lambda^\alpha\theta+u\cdot\nabla\theta=f $, 
  then for all $t\geq0$, we have 
  \begin{equation*}
    {\|\theta(t)\|}_{L^p}\leq {\|\theta(0)\|}_{L^p}+\int_{0}^{t}{\|f(\tau)\|}~{\rm d}\tau. 
  \end{equation*}
\end{lem}

The following lemma is a corollary of Theorem 1.1 and Lemma 3.3 in \cite{Ch_Mi_Zh_2007}. 

\begin{lem}\label{0325-6}
  Let $j\in\mathbb{N}$, $0< \alpha\leq 2$ and $2\leq p <\infty$. 
  Let $\theta,u$ and $f$ be smooth functions on $(0,\infty)\times\mathbb{R}^2$ 
  such that $\nabla\cdot u=0$. 
  If $\theta$ satisfies $\partial_t(\phi_j*\theta)+\Lambda^\alpha(\phi_j*\theta)+u\cdot\nabla(\phi_j*\theta)=f $, 
  then there exists constant $c_p>0$ such that for all $t\geq0$, we have 
  \begin{equation*}
    \partial_t{\|\phi_j*\theta(t)\|}_{L^p}+c_p2^{\alpha j}{\|\phi_j*\theta(t)\|}_{L^p}\leq {\|f(t)\|}_{L^p}. 
  \end{equation*}
\end{lem}

\begin{lem}[\cite{I_2020}]\label{0326-7}
  Let $j\in\mathbb{N}$ and $0< \alpha\leq 2$. 
  Let $\theta,u$ and $f$ be smooth functions on $(0,\infty)\times\mathbb{R}^d$ 
  such that $\theta, \partial_t\theta, \partial_t^2\theta, u\in L^\infty((0,\infty)\times\mathbb{R}^d)$. 
  If $\theta$ satisfies $\partial_t(\phi_j*\theta)+\Lambda^\alpha(\phi_j*\theta)+u\cdot\nabla(\phi_j*\theta)=f $, 
  then there exists constant $c>0$ such that for almost every $t\geq0$, we have 
  \begin{equation*}
    \partial_t{\|\phi_j*\theta(t)\|}_{L^\infty}+c2^{\alpha j}{\|\phi_j*\theta(t)\|}_{L^\infty}\leq {\|f(t)\|}_{L^\infty}. 
  \end{equation*}
\end{lem}

Thanks to Lemma \ref{0327-1}, 
we obtain convergence of linear solution. 

\begin{lem}\label{0326-3}
  Let $s\in\mathbb{R}$, $0<\alpha\leq2$, $1\leq p\leq\infty$ and $1\leq q,r<\infty$. 
  Then for any $\theta_0\in B_{p,q}^s$, we have 
  \begin{equation*}
    \lim_{T\to 0}{\|e^{-t\Lambda^\alpha}\theta_0\|}_{\widetilde{L}^r(0,T;B_{p,q}^{s+\alpha/r})}=0.
  \end{equation*}
\end{lem}

\begin{lem}[\cite{I_2020}]\label{0326-10}
  Let $s\in\mathbb{R}$, $0<\alpha\leq2$, $1\leq p\leq\infty$ and $1\leq r<\infty$. 
  If $\theta_0\in B_{p,\infty}^s$ satisfies that 
  \begin{equation*}
    \lim_{j\to\infty}2^{sj}{\|\phi_j*\theta_0\|}_{L^p}=0, 
  \end{equation*}
  then we have 
  \begin{equation*}
    \lim_{T\to 0}{\|e^{-t\Lambda^\alpha}\theta_0\|}_{\widetilde{L}^r(0,T;B_{p,\infty}^{s+\alpha/r})}=0.
  \end{equation*}
\end{lem}

Also, next lemma implies convergence of the nonlinear team in a mild solution of \eqref{SQG}. 

\begin{lem}\label{0326-4}
  Let $0<\alpha\leq2$, $1\leq p,q\leq\infty$ and $\theta\in B_{p,q}^{1-\alpha+2/p}$. 
  If $\theta$ is a mild solution of \eqref{SQG} for $\theta_0$, 
  then we have 
  \begin{equation*}
    \lim_{T\to 0}{\|\theta-e^{-t\Lambda^\alpha}\theta_0\|}_{L^\infty(0,T;B_{p,q}^{1-\alpha+2/p})}=0. 
  \end{equation*}
\end{lem}

\section{Proof of Theorem}
Let $\theta^{(1)},\theta^{(2)}$ are mild solutions of \eqref{SQG} with same initial data 
and we set 
\begin{equation*}
  w:=\theta^{(1)}-\theta^{(2)}.
\end{equation*}
By Lemma \ref{0301-33}, for any $j\in\mathbb{N}$, $w_j=\phi_j*w$ satisfies 
\begin{equation*}
   \begin{cases}
     \partial_t w_j + \Lambda^\alpha w_j + 
     \nabla\phi_j*\big((\nabla ^\perp \Lambda^{-1}\theta^{(1)}) \theta^{(1)}-(\nabla ^\perp \Lambda^{-1}\theta^{(2)}) \theta^{(2)}\big) = 0,\\
     w(0,x) = 0. 
   \end{cases}
\end{equation*}
and we introduce $N(\cdot, \cdot)$ by 
\begin{align}\label{0301-31}
   (\nabla ^\perp \Lambda^{-1}\theta^{(1)}) \theta^{(1)}-(\nabla ^\perp \Lambda^{-1}\theta^{(2)}) \theta^{(2)} 
   =&\frac{1}{2}\sum_{i=1}^{2}\big((\nabla ^\perp \Lambda^{-1} w)\theta^{(i)}+(\nabla ^\perp \Lambda^{-1}\theta^{(i)}) w\big)\notag\\ 
   =&:\frac{1}{2}\sum_{i=1}^{2}N(w,\theta^{(i)}). 
\end{align}
Thus, we have 
\begin{equation}\label{0227-5}
   \partial_t w_j + \Lambda^\alpha w_j + \frac{1}{2}\sum_{i=1}^{2}\nabla\phi_j*N(w,\theta^{(i)}) = 0. 
\end{equation} 
Similariy, we get 
\begin{equation}\label{0214-1}
   \partial_t (\psi*w) + \Lambda^\alpha (\psi*w) + \frac{1}{2}\sum_{i=1}^{2}\nabla\psi*N(w,\theta^{(i)}) = 0.  
\end{equation}
We devide into four cases depending on $\alpha$, 
\begin{equation*}
  \frac{3}{2}<\alpha\leq 2,\ \alpha=1,\ 1<\alpha\leq \frac{3}{2} \text{ and }0<\alpha<1. 
\end{equation*}
For a estimate of $\psi*w$, we prove only when the first case, $3/2<\alpha\leq 2$, 
since the frequency of $\psi*w$ exists near the origin and it is easier to estimate than $\phi_j*w$ ($j\in\mathbb{N}$). 

\subsection{Case $3/2<\alpha\leq 2$} 
Suppose $\theta^{(1)},\theta^{(2)}\in C([0,T];B_{p,p/(p-2)}^{-1/2})$. 
We show that 
\begin{equation*}
   {\|w\|}_{L^{p/2}(0,T_0;B_{p,p/2}^{-1/2})}=0\text{ for some small }T_0>0. 
\end{equation*}
Note that $p/2$ is the H\"older conjugate exponent of $p/(p-2)$ 
and since $p=4/(2\alpha-3)$, we have 
\begin{equation*}
   1< \frac{p}{p-2}\leq \frac{p}{2}< \infty. 
\end{equation*} 
By applying Lamma \ref{0325-1} to \eqref{0214-1}, $\psi*w$ satisfies 
\begin{equation*}
   {\|\psi*w(t)\|}_{L^p}\leq \frac{1}{2}\sum_{i=1}^{2}\int_{0}^{t}{\|\nabla\psi*N(w,\theta^{(i)})(\tau)\|}_{L^p}~{\rm d}\tau.
\end{equation*} 
Taking $L^{p/2}$ norm over the interval $[0,T]$, we get 
\begin{equation*}
   {\|\psi*w\|}_{L^{p/2}(0,T;L^p)}\leq \frac{1}{2}T\sum_{i=1}^{2}{\|\nabla\psi*N(w,\theta^{(i)})\|}_{L^{p/2}(0,T;L^p)}. 
\end{equation*}
By Bony's decomposition (see \cite{Bo_1981}) and Lemma \ref{0325-2}, we obtain 
\begin{align*}
   &{\|\nabla\psi*N(w,\theta^{(i)})\|}_{L^{p/2}(0,T;L^p)}\notag\\
   \leq& {\|\nabla \cdot N(w,\theta^{(i)})\|}_{L^{p/2}(0,T;B_{p,{p/2}}^{-1/2-\alpha})}\notag\\
   \leq& C \Bigg({\Bigg\|\sum_{l\geq 2}N(\psi*w, \theta^{(i)}_l)
   +\sum_{k\leq l-2}N(w_k, \theta^{(i)}_l)\Bigg\|}_{L^{p/2}(0,T;B_{p,p/2}^{-1/2-\alpha+1})}\notag\\
   &+{\Bigg\|\sum_{k\geq 2}N(w_k, \psi*\theta^{(i)})
   +\sum_{l \leq k-2}N(w_k, \theta^{(i)}_l)\Bigg\|}_{L^{p/2}(0,T;B_{p,p/2}^{-1/2-\alpha+1})}\notag\\
   &+{\|N(\psi*w, \psi*\theta^{(i)})
   +N(\psi*w, \phi_1*\theta^{(i)})
   +N(\phi_1*w, \psi*\theta^{(i)})\|}_{L^{p/2}(0,T;B_{p,p/2}^{-1/2-\alpha+1})}\notag\\
   &+{\Bigg\|\sum_{|k-l| \leq 1}\nabla\cdot N(w_k, \theta^{(i)}_l)\Bigg\|}_{L^{p/2}(0,T;B_{p,p/2}^{-1/2-\alpha})}\Bigg)\notag\\
   =:&C\big(({\rm I}\text{-}1)+({\rm I}\text{-}2)+({\rm I}\text{-}3a)+({\rm I}\text{-}3b)\big). 
\end{align*}
Using Lemma \ref{0325-3}, Lemma \ref{0325-4} and boundness of Riesz transform in $L^p$ ($1<p<\infty$), we have 
\begin{align}\label{0227-1}
   ({\rm I}\text{-}1)
   \leq& C\left({\|\nabla ^\perp \Lambda^{-1} w\|}_{L^{p/2}(0,T;B^{1-\alpha}_{\infty,\infty})}{\|\theta^{(i)}\|}_{L^\infty(0,T;B^{-1/2}_{p,p/2})}\right.\notag\\
   &\left.+{\|\nabla ^\perp \Lambda^{-1}\theta^{(i)}\|}_{L^\infty(0,T;B^{-1/2}_{p,p/2})}{\|w\|}_{L^{p/2}(0,T;B^{1-\alpha}_{\infty,\infty})}\right)\notag\\ 
   \leq& C{\|\theta^{(i)}\|}_{L^\infty(0,T;B^{-1/2}_{p,p/(p-2)})}{\|w\|}_{L^{p/2}(0,T;B^{-1/2}_{p,p/2})}. 
\end{align}
Similariy, we get 
\begin{equation}\label{0227-2}
   ({\rm I}\text{-}2)
   \leq C{\|\theta^{(i)}\|}_{L^\infty(0,T;B^{-1/2}_{p,p/(p-2)})}{\|w\|}_{L^{p/2}(0,T;B^{-1/2}_{p,p/2})}. 
\end{equation}
Since $\psi$ and $\phi_1$ has the support around the origin, 
we obtain 
\begin{equation}\label{0227-3}
   ({\rm I}\text{-}3a)
   \leq C{\|\theta^{(i)}\|}_{L^\infty(0,T;B^{-{1/2}}_{p,{p/(p-2)}})}{\|w\|}_{L^{p/2}(0,T;B^{-{1/2}}_{p,{p/2}})}. 
\end{equation}
By Lemma \ref{0325-4} and Lemma \ref{0325-5} \eqref{0326-2}, we have 
\begin{align}\label{0227-4}
   ({\rm I}\text{-}3b)
   &\leq C{\Bigg\|\sum_{|k-l| \leq 1}\nabla\cdot N(w_k, \theta^{(i)}_l)\Bigg\|}_{L^{p/2}(0,T;B_{p/2,p/2}^{-2})}\notag\\ 
   &\leq C{\|\theta^{(i)}\|}_{L^\infty(0,T;B^{-{1/2}}_{p,{p/(p-2)}})}{\|w\|}_{L^{p/2}(0,T;B^{-{1/2}}_{p,{p/2}})}. 
\end{align} 
From \eqref{0227-1}-\eqref{0227-4}, we obtain 
\begin{equation}\label{0301-5}
   {\|\psi*w\|}_{L^{p/2}(0,T;L^p)}\leq CT\sum_{i=1}^{2}{\|\theta^{(i)}\|}_{L^\infty(0,T;B^{-{1/2}}_{p,{p/(p-2)}})}{\|w\|}_{L^{p/2}(0,T;B^{-{1/2}}_{p,{p/2}})}. 
\end{equation}

Next, we estimate $w_j$ ($j\in\mathbb{N}$). 
We devide $\theta^{(i)}$ into the linear and the nonlinear terms as follows: 
\begin{equation*}
   \theta^{(i)}=e^{-t\Lambda^\alpha}\theta_0+(\theta^{(i)}-e^{-t\Lambda^\alpha}\theta_0)
   =:e^{-t\Lambda^\alpha}\theta_0+\widetilde{\theta}^{(i)}.
\end{equation*} 
Then \eqref{0227-5} is rewritten as 
\begin{align}\label{0301-12}
  &\partial_t w_j + \Lambda^\alpha w_j + \big(\nabla ^\perp \Lambda^{-1}(e^{-t\Lambda^\alpha}\theta_0)\big)\cdot \nabla w_j\notag\\
  =&-[\nabla\phi_j*, \nabla ^\perp \Lambda^{-1}(e^{-t\Lambda^\alpha}\theta_0)]w
  -\nabla\phi_j*\big((\nabla ^\perp \Lambda^{-1} w) (e^{-t\Lambda^\alpha}\theta_0)\big)\notag\\
  &-\frac{1}{2}\sum_{i=1}^{2}\nabla\phi_j*N(w,\widetilde{\theta}^{(i)}). 
\end{align}
Using Lemma \ref{0325-6}, we get 
\begin{align*}
  \partial_t{\|w_j\|}_{L^p}+c2^{\alpha j}{\|w_j\|}_{L^p}
  \leq{\|A_j(e^{-t\Lambda^\alpha}\theta_0,w)\|}_{L^p}
  +\frac{1}{2}\sum_{i=1}^{2}{\|\nabla\phi_j*N(w,\widetilde{\theta}^{(i)})\|}_{L^p}, 
\end{align*}
where 
\begin{equation*}
   A_j(f,g):=[\nabla\phi_j*, \nabla ^\perp \Lambda^{-1}f]g
   +\nabla\phi_j*\big((\nabla ^\perp \Lambda^{-1} g)f\big). 
\end{equation*}
By multiplying $e^{c2^{\alpha j}t}$ and integrating with respect to $t$, we have 
\begin{align}\label{0317-1}
  {\|w_j(t)\|}_{L^p}\leq&\int_{0}^{t}e^{-c2^{\alpha j}(t-\tau)}{\|A_j(e^{-t\Lambda^\alpha}\theta_0,w)(\tau)\|}_{L^p}~{\rm d}\tau\notag\\
  &+\frac{1}{2}\sum_{i=1}^{2}\int_{0}^{t}e^{-c2^{\alpha j}(t-\tau)}{\|\nabla\phi_j*N(w,\widetilde{\theta}^{(i)})(\tau)\|}_{L^p}~{\rm d}\tau. 
\end{align}
Taking $L^{p/2}$ norm over the interval $[0,T]$ and using Young's inequality, we obtain 
\begin{align}\label{0301-1}
  &{\|w_j\|}_{L^{p/2}(0,T;L^p)}\notag\\ 
  \leq&C\Bigg(2^{-2\alpha j/p}{\|A_j(e^{-t\Lambda^\alpha}\theta_0,w)\|}_{L^1(0,T;L^p)}
  +\sum_{i=1}^{2}2^{-\alpha j}{\|\nabla\phi_j*N(w,\widetilde{\theta}^{(i)})\|}_{L^{p/2}(0,T;L^p)}\Bigg). 
\end{align} 
Multiplying by $2^{-j/2}$ on both sides of \eqref{0301-1} and taking $l^{p/2}$ norm for $j$, we have 
\begin{align*}
  &{\Big\|\big\{2^{-j/2}{\|w_j\|}_{L^{p/2}(0,T;L^p)}\big\}_{j\in\mathbb{N}}\Big\|}_{l^{p/2}}\notag\\
  \leq &C\Bigg({\Big\|\big\{2^{(-1/2-2\alpha/p)j}{\|A_j(e^{-t\Lambda^\alpha}\theta_0,w)\|}_{L^1(0,T;L^p)}\big\}_{j\in\mathbb{N}}\Big\|}_{l^{p/2}}\notag\\
  &+\sum_{i=1}^{2}{\Bigg\|{\Big\|\big\{2^{(-1/2-\alpha)j}{\|\nabla\phi_j*N(w,\widetilde{\theta}^{(i)})\|}_{L^p}\big\}_{j\in\mathbb{N}}\Big\|}_{l^{p/2}}\Bigg\|}_{L^{p/2}(0,T)}\Bigg)\notag\\ 
  =:&C({\rm I}\hspace{-1.2pt}{\rm I}+{\rm I}\hspace{-1.2pt}{\rm I}\hspace{-1.2pt}{\rm I}). 
\end{align*}
If $3/2<\alpha<2$, then 
\begin{equation*}
   -\frac{1}{2}-\frac{2}{p}\alpha+\frac{2}{p}+1
   =\frac{1}{2}-\alpha^2+\frac{3}{2}\alpha+\alpha-\frac{3}{2}
   =-(\alpha-2)(\alpha-1)>0. 
\end{equation*}
Thus, by the commutator estimate (Lemma \ref{0325-7}) and the product estimate (Lemma \ref{0325-8}), 
we get 
\begin{align}\label{0301-2}
  {\rm I}\hspace{-1.2pt}{\rm I}
  \leq& {\Big\|\big\{2^{(-1/2-2\alpha/p)j}{\|[\phi_j, \nabla ^\perp \Lambda^{-1}(e^{-t\Lambda^\alpha}\theta_0)]\cdot \nabla w\|}_{L^1(0,T;L^p)}\big\}_{j\in\mathbb{N}}\Big\|}_{l^{p/2}}\notag\\ 
  &+{\Big\|\big\{2^{(-1/2-2\alpha/p)j}{\big\|\phi_j*\big((\nabla ^\perp \Lambda^{-1} w)\cdot \nabla (e^{-t\Lambda^\alpha}\theta_0)\big)\big\|}_{L^1(0,T;L^p)}\big\}_{j\in\mathbb{N}}\Big\|}_{l^{p/2}}\notag\\
  \leq&C{\|e^{-t\Lambda^\alpha}\theta_0\|}_{L^{p/(p-2)}(0,T;B_{p,p/(p-2)}^{-2\alpha/p+2/p+1})}
  {\|w\|}_{L^{p/2}(0,T;B_{p,p/2}^{-1/2})}. 
\end{align}
If $\alpha=2$, then $p=4$ and we introduce the following: 
\begin{equation*}
  ({\rm I}\hspace{-1.2pt}{\rm I}\text{-}1)
  :={\Bigg\|\Bigg\{2^{-3/2j}{\Bigg\|\sum_{k\geq 2}A_j(\psi*(e^{t\Delta}\theta_0),w_k)
  +\sum_{l \leq k-2}A_j((e^{t\Delta}\theta_0)_l,w_k)\Bigg\|}_{L^1(0,T;L^4)}\Bigg\}_{j\in\mathbb{N}}\Bigg\|}_{l^2}, 
\end{equation*}
\begin{equation*}
  ({\rm I}\hspace{-1.2pt}{\rm I}\text{-}2)
  :={\Bigg\|\Bigg\{2^{-3/2j}{\Bigg\|\sum_{l\geq 2}A_j((e^{t\Delta}\theta_0)_l,\psi*w)
  +\sum_{k\leq l-2}A_j((e^{t\Delta}\theta_0)_l,w_k)\Bigg\|}_{L^1(0,T;L^4)}\Bigg\}_{j\in\mathbb{N}}\Bigg\|}_{l^2}, 
\end{equation*}
\begin{equation*}
  ({\rm I}\hspace{-1.2pt}{\rm I}\text{-}3a)
  :={\Big\|\big\{2^{-3/2j}{\|A_j(\psi*(e^{t\Delta}\theta_0),\psi*w)\|}_{L^1(0,T;L^4)}\big\}_{j\in\mathbb{N}}\Big\|}_{l^2}, 
\end{equation*}
\begin{equation*}
  ({\rm I}\hspace{-1.2pt}{\rm I}\text{-}3b)
  :={\Big\|\big\{2^{-3/2j}{\|A_j(\psi*(e^{t\Delta}\theta_0),\phi_1*w)\|}_{L^1(0,T;L^4)}\big\}_{j\in\mathbb{N}}\Big\|}_{l^2}, 
\end{equation*}
\begin{equation*}
  ({\rm I}\hspace{-1.2pt}{\rm I}\text{-}3c)
  :={\Big\|\big\{2^{-3/2j}{\|A_j(\phi_1*w(e^{t\Delta}\theta_0),\psi*w)\|}_{L^1(0,T;L^4)}\big\}_{j\in\mathbb{N}}\Big\|}_{l^2}, 
\end{equation*}
\begin{equation*}
  ({\rm I}\hspace{-1.2pt}{\rm I}\text{-}4)
  :={\Bigg\|\Bigg\{2^{-3/2j}{\Bigg\|\sum_{|k-l| \leq 1}(\nabla ^\perp \Lambda^{-1}(e^{t\Delta}\theta_0)_l) \cdot\nabla(w_j)_k\Bigg\|}_{L^1(0,T;L^4)}\Bigg\}_{j\in\mathbb{N}}\Bigg\|}_{l^2}, 
\end{equation*}
\begin{equation*}
  ({\rm I}\hspace{-1.2pt}{\rm I}\text{-}5)
  :={\Bigg\|\Bigg\{2^{-3/2j}{\Bigg\|\sum_{|k-l| \leq 1}\nabla\phi_j*N(w_k,(e^{t\Delta}\theta_0)_l)\Bigg\|}_{L^1(0,T;L^4)}\Bigg\}_{j\in\mathbb{N}}\Bigg\|}_{l^2}. 
\end{equation*}
Then we can write 
\begin{align*}
  {\rm I}\hspace{-1.2pt}{\rm I}
  \leq({\rm I}\hspace{-1.2pt}{\rm I}\text{-}1)
  +({\rm I}\hspace{-1.2pt}{\rm I}\text{-}2)
  +({\rm I}\hspace{-1.2pt}{\rm I}\text{-}3a)+({\rm I}\hspace{-1.2pt}{\rm I}\text{-}3b)+({\rm I}\hspace{-1.2pt}{\rm I}\text{-}3c)
  +({\rm I}\hspace{-1.2pt}{\rm I}\text{-}4)
  +({\rm I}\hspace{-1.2pt}{\rm I}\text{-}5). 
\end{align*}
Applying Lemma \ref{0325-3} and mean-value theorem, we have 
\begin{equation*}
   ({\rm I}\hspace{-1.2pt}{\rm I}\text{-}1),({\rm I}\hspace{-1.2pt}{\rm I}\text{-}2)
   \leq C{\|e^{t\Delta}\theta_0\|}_{L^2(0,T;B_{4,2}^{1/2})}{\|w\|}_{L^2(0,T;B_{4,2}^{-1/2})}.
\end{equation*}
By a similar estimate to \eqref{0227-3}, we obtain 
\begin{equation*}
   ({\rm I}\hspace{-1.2pt}{\rm I}\text{-}3a)\text{-}({\rm I}\hspace{-1.2pt}{\rm I}\text{-}3c)
   \leq C{\|e^{t\Delta}\theta_0\|}_{L^2(0,T;B_{4,2}^{1/2})}{\|w\|}_{L^2(0,T;B_{4,2}^{-1/2})}. 
\end{equation*}
Since $f_{j,k}:=\phi_j*\phi_k*f=0$ (if $|j-k|\geq 2$), we have 
\begin{align*}
  ({\rm I}\hspace{-1.2pt}{\rm I}\text{-}4)&
  \leq \sum_{j\in\mathbb{N}}2^{-j/2}\sum_{k\in\mathbb{N}, |j-k|\leq 1}\ \sum_{l\in\mathbb{N}, |k-l|\leq 1}{\|\nabla ^\perp \Lambda^{-1}(e^{t\Delta}\theta_0)_l w_{j,k}\|}_{L^1(0,T;L^4)}\\
  &\leq C\sum_{j\in\mathbb{N}}\sum_{k\in\mathbb{N}, |j-k|\leq 1}2^{-(j-k)/2}\sum_{l\in\mathbb{N}, |k-l|\leq 1}{\|\nabla ^\perp \Lambda^{-1}(e^{t\Delta}\theta_0)_l\|}_{L^2(0,T;L^\infty)}2^{-k/2} {\|w_k\|}_{L^2(0,T;L^2)}\\
  &\leq C{\|e^{t\Delta}\theta_0\|}_{L^{2}(0,T;B_{4,2}^{1/2})}{\|w\|}_{L^2(0,T;B_{4,2}^{-1/2})}. 
\end{align*}
Using Lemma \ref{0325-5} \eqref{0326-1}, we get 
\begin{equation*}
   ({\rm I}\hspace{-1.2pt}{\rm I}\text{-}5)
  \leq C{\|e^{t\Delta}\theta_0\|}_{L^2(0,T;B_{4,2}^{1/2})}{\|w\|}_{L^2(0,T;B_{4,2}^{-1/2})}, 
\end{equation*}
and 
\begin{equation}\label{0301-3}
  {\rm I}\hspace{-1.2pt}{\rm I}
  \leq C{\|e^{t\Delta}\theta_0\|}_{L^2(0,T;B_{4,2}^{1/2})}{\|w\|}_{L^2(0,T;B_{4,2}^{-1/2})}. 
\end{equation}
Using a similar estimate to $\nabla\psi*N(w,\theta^{(i)})$, we have 
\begin{align}\label{0301-4}
  {\rm I}\hspace{-1.2pt}{\rm I}\hspace{-1.2pt}{\rm I}
  \leq &C\sum_{i=1}^{2}{\|\widetilde{\theta}^{(i)}\|}_{L^\infty(0,T;B^{-1/2}_{p,p/(p-2)})}{\|w\|}_{L^{p/2}(0,T;B^{-1/2}_{p,p/2})}. 
\end{align}
From \eqref{0301-2}-\eqref{0301-4}, we obtain 
\begin{align}\label{0301-6}
  &{\Big\|{\big\|\big\{2^{-j/2}{\|w_j\|}_{L^p}\big\}_{j\in\mathbb{N}}\big\|}_{l^{p/2}}\Big\|}_{L^{p/2}(0,T)}\notag\\
  \leq& C\left({\|e^{-t\Lambda^\alpha}\theta_0\|}_{L^{p/(p-2)}(0,T;B_{p,p/(p-2)}^{-2\alpha/p+2/p+1})}+\sum_{i=1}^{2}{\|\widetilde{\theta}^{(i)}\|}_{L^\infty(0,T;B^{-1/2}_{p,p/(p-2)})}\right)
  {\|w\|}_{L^{p/2}(0,T;B^{-1/2}_{p,p/2})} 
\end{align}
Summing \eqref{0301-5} and \eqref{0301-6}, we get 
\begin{align*}
  {\|w\|}_{L^{p/2}(0,T;B^{-1/2}_{p,p/2})}
  \leq &C\Bigg({\|e^{-t\Lambda^\alpha}\theta_0\|}_{L^{p/(p-2)}(0,T;B_{p,p/(p-2)}^{-2\alpha/p+2/p+1})}\\ 
  &+\sum_{i=1}^{2}\Big(T{\|\theta^{(i)}\|}_{L^\infty(0,T;B^{-1/2}_{p,p/(p-2)})}
  +{\|\widetilde{\theta}^{(i)}\|}_{L^\infty(0,T;B^{-1/2}_{p,p/(p-2)})}\Big)\Bigg)\notag\\ 
  &\times{\|w\|}_{L^{p/2}(0,T;B^{-1/2}_{p,p/2})}. 
\end{align*}
Now, we have 
\begin{equation*}
   -\frac{2}{p}\alpha+\frac{2}{p}+1
   =-\frac{2}{p}\alpha+\alpha-\frac{3}{2}+1
   =-\frac{1}{2}+\frac{p-2}{p}\alpha. 
\end{equation*}
Therefore, since Lemma \ref{0326-3} and Lemma \ref{0326-4}, 
for any $\varepsilon>0$, there exists $0<T_0\ll 1$ such that 
\begin{equation*}
   {\|w\|}_{L^{p/2}(0,{T_0};B^{-1/2}_{p,p/2})}
   \leq C\varepsilon{\|w\|}_{L^{p/2}(0,{T_0};B^{-1/2}_{p,p/2})}. 
\end{equation*}
By taking $\varepsilon$ sufficiently small, we obtain 
\begin{equation*}
   {\|w\|}_{L^{p/2}(0,{T_0};B^{-1/2}_{p,p/2})}
   \leq 0. 
\end{equation*}
We conclude that 
\begin{equation*}
   \theta^{(1)}=\theta^{(2)}\ \text{in}\ B^{-1/2}_{p,p/2}\ \text{for}\ t\in[0,T_0].
\end{equation*}
We prove that ${\|w\|}_{B^{-1/2}_{p,p/2}}=0$\ in the entire interval $[0,T]$, by contradiction argument.
Set
$$\tau^*:=\sup\{\tau\in[0,T)\ |\ {\|\theta^{(1)}(t,\cdot)-\theta^{(2)}(t,\cdot)\|}_{B^{-1/2}_{p,p/2}}=0\ \text{for}\ t\in[0,\tau]\}.$$
Assume that $\tau^*<T$. 
From the continuity in time of $\theta^{(1)}$ and $\theta^{(2)}$, we have $\theta^{(1)}(\tau^*)=\theta^{(2)}(\tau^*).$
By the same argument as local uniqueness , there exists $\delta'>0$ such that $\theta^{(1)}(s+\tau^*)=\theta^{(2)}(s+\tau^*)$ for $0\leq s\leq\delta'$.
This is a contradiction with the definition of $\tau^*$.
Therefore $\tau^*=T$. 

\subsection{Case $\alpha=1$} 
In this case, $\theta^{(1)},\theta^{(2)}\in C([0,T];B_{\infty,\infty}^0)$. 
Let $s\in\mathbb{R}$ with $0<s<1$. 
We consider ${\|w\|}_{L^\infty(0,T;B_{\infty,\infty}^{-s})}$. 
We devide $\theta^{(i)}$ as follows: 
\begin{equation*}
  \theta^{(i)}=S_2\theta^{(i)}+\widetilde{S}_2\theta^{(i)}, 
\end{equation*}
where $S_2$ and $\widetilde{S}_2$ are defined by \eqref{0308-3}. 
Note that 
\begin{equation*}
   {\|S_2f\|}_{B_{\infty,\infty}^{-s}}\leq C{\|f\|}_{B_{\infty,\infty}^{-s}} 
   \text{ and }\ {\|\widetilde{S}_2f\|}_{B_{\infty,\infty}^{-s}}\leq C{\|f\|}_{B_{\infty,\infty}^{-s}}. 
\end{equation*}
For a estimate of $\psi*w$, 
using Lemma \ref{0325-7} and Lemma \ref{0325-8} for $\widetilde{S}_2\theta^{(i)}$ and $\widetilde{S}_2\theta^{(i)}$ respectively, 
we obtain 
\begin{align}\label{0415-6}
  &{\|\psi*w\|}_{L^\infty(0,T;L^\infty)}\notag\\
  \leq& C\sum_{i=1}^{2}T{\|\theta^{(i)}\|}_{L^\infty(0,T;B^{0}_{\infty,\infty})}
  \big({\|\nabla^\perp\Lambda^{-1}(\psi*w)\|}_{L^\infty(0,T;L^\infty)}+{\|w\|}_{L^\infty(0,T;B_{\infty, \infty}^{-s})}\big). 
\end{align}

We consider $w_j$ ($j\in\mathbb{N}$). 
For each $j$, we devide the nonlinear term $\widetilde{\theta}^{(i)}$ as follows: 
\begin{equation*}
  \widetilde{\theta}^{(i)}
  =S_{j+2}\widetilde{\theta}^{(i)}+\widetilde{S}_{j+2}\widetilde{\theta}^{(i)}. 
\end{equation*}
Then \eqref{0301-12} is rewritten as 
\begin{align*}
  &\partial_t w_j + \Lambda w_j + \Big(\nabla ^\perp \Lambda^{-1}(e^{-t\Lambda}\theta_0)+\frac{1}{2}\sum_{i=1}^{2}\nabla ^\perp \Lambda^{-1}(S_{j+2}\widetilde{\theta}^{(i)})\Big)\cdot \nabla w_j\\
  =&-[\nabla\phi_j*, \nabla ^\perp \Lambda^{-1}(e^{-t\Lambda}\theta_0)]w
  -\nabla\phi_j*\big((\nabla ^\perp \Lambda^{-1} w)e^{-t\Lambda}\theta_0\big)\\
  &-\frac{1}{2}\sum_{i=1}^{2}\Big([\nabla\phi_j*, \nabla ^\perp \Lambda^{-1}(S_{j+2}\widetilde{\theta}^{(i)})]w
  +\nabla\phi_j*\big((\nabla ^\perp \Lambda^{-1} w)S_{j+2}\widetilde{\theta}^{(i)}\big)\\
  &+\nabla\phi_j*N\big(w,\widetilde{S}_{j+2}\widetilde{\theta}^{(i)}\big)\Big). 
\end{align*}
Using the maximum principle (Lemma \ref{0326-7}) , we have 
\begin{align*}
  &\partial_t{\|w_j\|}_{L^\infty}+c2^{j}{\|w_j\|}_{L^\infty}\notag\\ 
  \leq&{\|[\nabla\phi_j*, \nabla ^\perp \Lambda^{-1}(e^{-t\Lambda}\theta_0)]w\|}_{L^\infty} 
  +{\big\|\nabla\phi_j*\big((\nabla ^\perp \Lambda^{-1} w)e^{-t\Lambda}\theta_0\big)\big\|}_{L^\infty}\notag\\
  &+\frac{1}{2}\sum_{i=1}^{2}\Big({\|[\nabla\phi_j*, \nabla ^\perp \Lambda^{-1}(S_{j+2}\widetilde{\theta}^{(i)})]w\|}_{L^\infty} 
  +{\big\|\nabla\phi_j*\big((\nabla ^\perp \Lambda^{-1} w)S_{j+2}\widetilde{\theta}^{(i)}\big)\big\|}_{L^\infty}\notag\\ 
  &+{\big\|\nabla\phi_j*N\big(w,\widetilde{S}_{j+2}\widetilde{\theta}^{(i)}\big)\big\|}_{L^\infty}\Big).  
\end{align*}
Let $r\in\mathbb{R}$ with $1<r<1/s$. 
Using Gronwall's inequality, taking $L^\infty$ norm over the interval $[0,T]$, 
multiplying $2^{-sj}$ and taking $l^\infty$ norm for $j$, we obtain 
\begin{align}\label{0301-29}
  &{\Big\|\big\{2^{-sj}{\|w_j\|}_{L^{\infty}(0,T;L^\infty)}\big\}_{j\in\mathbb{N}}\Big\|}_{l^\infty}\notag\\
  \leq &C\Bigg({\Big\|\big\{2^{(-s-1+1/r)j}{\|[\nabla\phi_j*, \nabla ^\perp \Lambda^{-1}(e^{-t\Lambda}\theta_0)]w\|}_{L^r(0,T;L^\infty)}\big\}_{j\in\mathbb{N}}\Big\|}_{l^\infty}\notag\\
  &+{\Big\|\big\{2^{(-s-1+1/r)j}{\|\nabla\phi_j*\big((\nabla ^\perp \Lambda^{-1} w)e^{-t\Lambda}\theta_0\big)\|}_{L^r(0,T;L^\infty)}\big\}_{j\in\mathbb{N}}\Big\|}_{l^\infty}\notag\\
  &+\sum_{i=1}^{2}\Bigg({\Big\|\big\{2^{(-s-1)j}{\|[\nabla\phi_j*, \nabla ^\perp \Lambda^{-1}(S_{j+2}\widetilde{\theta}^{(i)})]w\|}_{L^\infty(0,T;L^\infty)}\big\}_{j\in\mathbb{N}}\Big\|}_{l^\infty}\notag\\
  &+{\Big\|\big\{2^{(-s-1)j}{\big\|\nabla\phi_j*\big((\nabla ^\perp \Lambda^{-1} w)S_{j+2}\widetilde{\theta}^{(i)}\big)\big\|}_{L^\infty(0,T;L^\infty)}\big\}_{j\in\mathbb{N}}\Big\|}_{l^\infty}\notag\\
  &+{\Bigg\|{\Big\|\big\{2^{(-s-1)j}{\big\|\nabla\phi_j*N\big(w,\widetilde{S}_{j+2}\widetilde{\theta}^{(i)}\big)\big\|}_{L^\infty}\big\}_{j\in\mathbb{N}}\Big\|}_{l^\infty}\Bigg\|}_{L^{\infty}(0,T)}\Bigg)\Bigg)\notag\\ 
  =:&C({\rm I}\hspace{-1.2pt}{\rm V}
  +{\rm V}
  +{\rm V}\hspace{-1.2pt}{\rm I}
  +{\rm V}\hspace{-1.2pt}{\rm I}\hspace{-1.2pt}{\rm I}
  +{\rm V}\hspace{-1.2pt}{\rm I}\hspace{-1.2pt}{\rm I}\hspace{-1.2pt}{\rm I}). 
\end{align}
Since $0<s<1$ and $1<r<1/s$, we have 
\begin{equation*}
   -s-1+\frac{1}{r}+1>0. 
\end{equation*}
Thus, using Lemma \ref{0325-7} and Lemma \ref{0325-8}, we get 
\begin{equation}\label{0415-1}
  {\rm I}\hspace{-1.2pt}{\rm V} 
  \leq C {\|e^{-t\Lambda}\theta_0\|}_{\widetilde{L}^r(0,T;B_{\infty,\infty}^{1/r})} 
  {\|w\|}_{L^\infty(0,T;B_{\infty,\infty}^{-s})} 
\end{equation}
and 
\begin{equation}\label{0415-2}
  {\rm V} 
  \leq C {\|e^{-t\Lambda}\theta_0\|}_{\widetilde{L}^r(0,T;B_{\infty,\infty}^{1/r})}
  \big({\|\nabla^\perp\Lambda^{-1}(\psi*w)\|}_{L^\infty(0,T;L^\infty)}+{\|w\|}_{L^\infty(0,T;B_{\infty,\infty}^{-s})}\big). 
\end{equation}
We estimate V\hspace{-1.2pt}I. 
Since $\phi_j*\big(\big(S_{j+2}f\big)g_k\big)=0$ (if $k>j+1$), we get 
\begin{align*}
  [\nabla\phi_j*, \nabla ^\perp \Lambda^{-1}(S_{j+2}\widetilde{\theta}^{(i)})]w
  =&[\nabla\phi_j*, \nabla ^\perp \Lambda^{-1}(S_{j+2}\widetilde{\theta}^{(i)})](\psi*w)\\
  &+\sum_{k\leq j+1}[\nabla\phi_j*, \nabla ^\perp \Lambda^{-1}(S_{j+2}\widetilde{\theta}^{(i)})]w_k. 
\end{align*}
Thus, due to mean value theorem, we have 
\begin{align*}
  {\rm V}\hspace{-1.2pt}{\rm I}
  \leq&\sum_{i=1}^{2}\sup_{j\in\mathbb{N}}2^{(-s-1)j}2^{-j}\Bigg({\|\nabla\nabla^\perp\Lambda^{-1}(S_{j+2}\widetilde{\theta}^{(i)})\|}_{L^\infty(0,T;L^\infty)}{\|\nabla(\psi*w)\|}_{L^\infty(0,T;L^\infty)}\notag\\
  &+\sum_{k\leq j+1}{\|\nabla\nabla^\perp\Lambda^{-1}(S_{j+2}\widetilde{\theta}^{(i)})\|}_{L^\infty(0,T;L^\infty)}
  {\|\nabla w_k\|}_{L^\infty(0,T;L^\infty)}\Bigg). 
\end{align*}
Using Young's inequality, we obtain 
\begin{equation}\label{0415-3}
  {\rm V}\hspace{-1.2pt}{\rm I} 
  \leq C {\|\widetilde{\theta}^{(i)}\|}_{L^\infty(0,T;B_{\infty,\infty}^0)}
  {\|w\|}_{L^\infty(0,T;B_{\infty,\infty}^{-s})} 
\end{equation}
Similariy, we get 
\begin{equation}\label{0415-4}
  {\rm V}\hspace{-1.2pt}{\rm I}\hspace{-1.2pt}{\rm I} 
  \leq C{\|\widetilde{\theta}^{(i)}\|}_{L^\infty(0,T;B_{\infty,\infty}^0)}
  \big({\|\nabla^\perp\Lambda^{-1}(\psi*w)\|}_{L^\infty(0,T;L^\infty)}+{\|w\|}_{L^\infty(0,T;B_{\infty,\infty}^{-s})}\big). 
\end{equation}
We consider {\rm V}\hspace{-1.2pt}{\rm I}\hspace{-1.2pt}{\rm I}\hspace{-1.2pt}{\rm I}. 
Since $\psi*(\widetilde{S}_{j+2}\widetilde{\theta}^{(i)})=\widetilde{S}_{j+2}\widetilde{\theta}_1^{(i)}=0$ for any $j\in\mathbb{N}$, 
by Bony's decomposition, we have 
\begin{align*}
  {\rm V}\hspace{-1.2pt}{\rm I}\hspace{-1.2pt}{\rm I}\hspace{-1.2pt}{\rm I}
  \leq&\sum_{i=1}^{2} \left({\left\|{\Bigg\|\Bigg\{2^{(-s-1)j}{\Bigg\|\sum_{l\geq 2}\nabla\phi_j*N\big(\psi*w,\widetilde{S}_{j+2}\widetilde{\theta}_l^{(i)}\big)\Bigg\|}_{L^\infty}\Bigg\}_{j\in\mathbb{N}}\Bigg\|}_{l^\infty}\right\|}_{L^{\infty}(0,T)}\right.\notag\\
  &+{\left\|{\Bigg\|\Bigg\{2^{(-s-1)j}{\Bigg\|\sum_{k\leq l-2}\nabla\phi_j*N\big(w_k,\widetilde{S}_{j+2}\widetilde{\theta}_l^{(i)}\big)\Bigg\|}_{L^\infty}\Bigg\}_{j\in\mathbb{N}}\Bigg\|}_{l^\infty}\right\|}_{L^{\infty}(0,T)}\notag\\
  &+{\left\|{\Bigg\|\Bigg\{2^{(-s-1)j}{\Bigg\|\sum_{l \leq k-2}\nabla\phi_j*N\big(w_k,\widetilde{S}_{j+2}\widetilde{\theta}_l^{(i)}\big)\Bigg\|}_{L^\infty}\Bigg\}_{j\in\mathbb{N}}\Bigg\|}_{l^\infty}\right\|}_{L^{\infty}(0,T)}\notag\\
  &+\left.{\left\|{\Bigg\|\Bigg\{2^{(-s-1)j}{\Bigg\|\sum_{|k-l| \leq 1}\nabla\phi_j*N\big(w_k,\widetilde{S}_{j+2}\widetilde{\theta}_l^{(i)}\big)\Bigg\|}_{L^\infty}\Bigg\}_{j\in\mathbb{N}}\Bigg\|}_{l^\infty}\right\|}_{L^{\infty}(0,T)}\right)\notag\\
  =:&({\rm V}\hspace{-1.2pt}{\rm I}\hspace{-1.2pt}{\rm I}\hspace{-1.2pt}{\rm I}\text{-}1a)
  +({\rm V}\hspace{-1.2pt}{\rm I}\hspace{-1.2pt}{\rm I}\hspace{-1.2pt}{\rm I}\text{-}1b)
  +({\rm V}\hspace{-1.2pt}{\rm I}\hspace{-1.2pt}{\rm I}\hspace{-1.2pt}{\rm I}\text{-}2)
  +({\rm V}\hspace{-1.2pt}{\rm I}\hspace{-1.2pt}{\rm I}\hspace{-1.2pt}{\rm I}\text{-}3). 
\end{align*}
Note that for any $j,k\in\mathbb{N}$, we have 
\begin{equation*}
   {\|\widetilde{S}_{j+2}f_k\|}_{L^\infty}\leq C{\|f_k\|}_{L^\infty}
   \text{ and }
   {\|\nabla^\perp\Lambda^{-1}(\widetilde{S}_{j+2}f_k)\|}_{L^\infty}\leq C{\|\widetilde{S}_{j+2}f_k\|}_{L^\infty}, 
\end{equation*}
where $C$ is independent of $j$. 
Thus, by a similar estimate to Lemma \ref{0325-3}, we obtain 
\begin{equation*}
  ({\rm V}\hspace{-1.2pt}{\rm I}\hspace{-1.2pt}{\rm I}\hspace{-1.2pt}{\rm I}\text{-}1a), ({\rm V}\hspace{-1.2pt}{\rm I}\hspace{-1.2pt}{\rm I}\hspace{-1.2pt}{\rm I}\text{-}1b)
  \leq C{\|\widetilde{\theta}^{(i)}\|}_{\widetilde{L}^\infty(0,T;B_{\infty,\infty}^{-s-1})}
  \big({\|\nabla^\perp\Lambda^{-1}(\psi*w)\|}_{L^\infty(0,T;L^\infty)}+{\|w\|}_{L^\infty(0,T;B_{\infty,\infty}^{-s-1})}\big). 
\end{equation*}
Thanks to a frequency is restricted to $\widetilde{S}_{j+2}$, 
we get 
\begin{align*}
  ({\rm V}\hspace{-1.2pt}{\rm I}\hspace{-1.2pt}{\rm I}\hspace{-1.2pt}{\rm I}\text{-}2)
  = &{\left\|\sup_{j\in\mathbb{N}}2^{(-s-1)j}{\Bigg\|\sum_{k\in\mathbb{N},|k-j|\leq 4}\ \sum_{j+2\leq l\leq k-2}\nabla\phi_j*N\big(w_k,\widetilde{S}_{j+2}\widetilde{\theta}_l^{(i)}\big)\Bigg\|}_{L^\infty}\right\|}_{L^{\infty}(0,T)}\\
  = &{\Bigg\|\sup_{j\in\mathbb{N}}2^{(-s-1)j}{\Big\|\nabla\phi_j*N\big(w_{j+4},\widetilde{S}_{j+2}\widetilde{\theta}_{j+2}^{(i)}\big)\Big\|}_{L^\infty}\Bigg\|}_{L^{\infty}(0,T)}. 
\end{align*}
Hence we obtain 
\begin{align*}
  ({\rm V}\hspace{-1.2pt}{\rm I}\hspace{-1.2pt}{\rm I}\hspace{-1.2pt}{\rm I}\text{-}2)
  \leq& C {\Bigg\|\sup_{j\in\mathbb{N}}2^{-sj}{\|w_{j+4}\|}_{L^\infty}{\big\|\widetilde{S}_{j+2}\widetilde{\theta}_{j+2}^{(i)}\big\|}_{L^\infty}\Bigg\|}_{L^{\infty}(0,T)}\\
  \leq& C {\|\widetilde{\theta}^{(i)}\|}_{L^\infty(0,T;B_{\infty,\infty}^0)}{\|w\|}_{L^\infty(0,T;B_{\infty,\infty}^{-s})}. 
\end{align*}
Using Lemma \ref{0325-5}, we have 
\begin{equation*}
  ({\rm V}\hspace{-1.2pt}{\rm I}\hspace{-1.2pt}{\rm I}\hspace{-1.2pt}{\rm I}\text{-}3)
  \leq C{\|\widetilde{\theta}^{(i)}\|}_{L^\infty(0,T;B_{\infty,\infty}^0)}
  {\|w\|}_{L^\infty(0,T;B_{\infty,\infty}^{-s})}, 
\end{equation*}
Thus, we get 
\begin{equation}\label{0415-5}
  {\rm V}\hspace{-1.2pt}{\rm I}\hspace{-1.2pt}{\rm I}\hspace{-1.2pt}{\rm I}
  \leq C{\|\widetilde{\theta}^{(i)}\|}_{L^\infty(0,T;B_{\infty,\infty}^{0})}
  \big({\|\nabla^\perp\Lambda^{-1}(\psi*w)\|}_{L^\infty(0,T;L^\infty)}+{\|w\|}_{L^\infty(0,T;B_{\infty,\infty}^{-s})}\big). 
\end{equation}
From \eqref{0415-1}-\eqref{0415-5}, we obtain 
\begin{align}\label{0415-7}
  {\Big\|\big\{2^{-sj}{\|w_j\|}_{L^{\infty}(0,T;L^\infty)}\big\}_{j\in\mathbb{N}}\Big\|}_{l^\infty}\notag 
  \leq& C\big({\|e^{-t\Lambda}\theta_0\|}_{\widetilde{L}^r(0,T;B_{\infty,\infty}^{1/r})}
  +{\|\widetilde{\theta}^{(i)}\|}_{L^\infty(0,T;B_{\infty,\infty}^{0})}\big)\notag\\ 
  &\times\big({\|\nabla^\perp\Lambda^{-1}(\psi*w)\|}_{L^\infty(0,T;L^\infty)}+{\|w\|}_{L^\infty(0,T;B_{\infty,\infty}^{-s})}\big). 
\end{align}

It remains to estimate ${\|\nabla^\perp\Lambda^{-1}(\psi*w)\|}_{L^\infty(0,T;L^\infty)}$. 
We perform the Riesz transform on both side of \eqref{0214-1}, 
then we have 
\begin{align*}
  &\partial_t \nabla^\perp\Lambda^{-1}(\psi*w)+ \Lambda \nabla^\perp\Lambda^{-1}(\psi*w) +\frac{1}{2}\sum_{i=1}^{2}(\nabla^\perp\Lambda^{-1}(S_2\theta^{(i)}))\cdot\nabla \nabla^\perp\Lambda^{-1}(\psi*w)\notag\\
  =&-\frac{1}{2}\sum_{i=1}^{2}\Big([\nabla\nabla^\perp\Lambda^{-1}\psi*,\nabla^\perp\Lambda^{-1}(S_2\theta^{(i)})]w 
  +\nabla\nabla^\perp\Lambda^{-1} \big(\psi*((\nabla^\perp\Lambda^{-1}w) (S_2\theta^{(i)}))\big)\notag\\
  &+\nabla\nabla^\perp\Lambda^{-1}\big(\psi*N\big(w,\widetilde{S}_2\theta^{(i)}\big)\big)\Big). 
\end{align*}
Using Lemma \eqref{0325-1}, we obtain 
\begin{align*}
  {\|\nabla^\perp\Lambda^{-1}(\psi*w)\|}_{L^\infty(0,T;L^\infty)}
  \leq& \frac{1}{2}T\sum_{i=1}^{2}\Big({\big\|[\nabla\nabla^\perp\Lambda^{-1}\psi*,\nabla^\perp\Lambda^{-1}(S_2\theta^{(i)})]w\big\|}_{L^\infty(0,T;L^\infty)}\\
  &+{\big\|\nabla\nabla^\perp\Lambda^{-1} \big(\psi*((\nabla^\perp\Lambda^{-1}w) (S_2\theta^{(i)}))\big)\big\|}_{L^\infty(0,T;L^\infty)}\\ 
  &+{\big\|\nabla\nabla^\perp\Lambda^{-1}\big(\psi*N(w,\widetilde{S}_2\theta^{(i)})\big)\big\|}_{L^\infty(0,T;L^\infty)}\Big). 
\end{align*}
Thanks to $S_2f_j=0$ (if $k>3$), by Lemma \ref{0326-9}, we have 
\begin{align*}
  &{\big\|[\nabla\nabla^\perp\Lambda^{-1}\psi*,\nabla^\perp\Lambda^{-1}(S_2\theta^{(i)})]w\big\|}_{L^\infty(0,T;L^\infty)}\\ 
  \leq& C \Bigg({\|\psi*(S_2\theta^{(i)})\|}_{L^\infty}{\|\psi*w\|}_{L^\infty} 
  +{\|\psi*(S_2\theta^{(i)})\|}_{L^\infty}{\|w_1\|}_{L^\infty}\\ 
  &+{\|S_2\theta_1^{(i)}\|}_{L^\infty}{\|\psi*w\|}_{L^\infty} 
  +\sum_{k=1}^{3}\sum_{l\in\mathbb{N},|k-l|\leq 6}{\|S_2\theta_k^{(i)}\|}_{L^\infty}{\|w_l\|}_{L^\infty}\Bigg)\\ 
  \leq& C{\|\theta^{(i)}\|}_{L^\infty(0,T;B_{\infty,\infty}^{0})}
  {\|w\|}_{L^\infty(0,T;B_{\infty,\infty}^{-s})}. 
\end{align*}
Using Lemma \ref{0326-5}, we get 
\begin{align*}
  &{\big\|\nabla\nabla^\perp\Lambda^{-1} \big(\psi*((\nabla^\perp\Lambda^{-1}w) (S_2\theta^{(i)}))\big)\big\|}_{L^\infty(0,T;L^\infty)}\\ 
  \leq& C{\big\|\psi*((\nabla^\perp\Lambda^{-1}w) (S_2\theta^{(i)}))\big\|}_{L^\infty(0,T;L^\infty)}\\ 
  \leq& C{\|\theta^{(i)}\|}_{L^\infty(0,T;B_{\infty,\infty}^{0})}
  \big({\|\nabla^\perp\Lambda^{-1}(\psi*w)\|}_{L^\infty(0,T;L^\infty)}+{\|w\|}_{L^\infty(0,T;B_{\infty,\infty}^{-s})}\big), 
\end{align*}
and 
\begin{align*}
  &{\big\|\nabla\nabla^\perp\Lambda^{-1}\big(\psi*N(w,\widetilde{S}_2\theta^{(i)})\big)\big\|}_{L^\infty(0,T;L^\infty)}\\ 
  \leq& C{\|\theta^{(i)}\|}_{L^\infty(0,T;B_{\infty,\infty}^{0})}
  \big({\|\nabla^\perp\Lambda^{-1}(\psi*w)\|}_{L^\infty(0,T;L^\infty)}+{\|w\|}_{L^\infty(0,T;B_{\infty,\infty}^{-s})}\big). 
\end{align*}
Therefore, we have 
\begin{align}\label{0301-27}
  &{\|\nabla^\perp\Lambda^{-1}(\psi*w)\|}_{L^\infty(0,T;L^\infty)}\notag\\
  \leq& CT\sum_{i=1}^{2}{\|\theta^{(i)}\|}_{L^\infty(0,T;B_{\infty,\infty}^{0})}
  \big({\|\nabla^\perp\Lambda^{-1}(\psi*w)\|}_{L^\infty(0,T;L^\infty)}+{\|w\|}_{L^\infty(0,T;B_{\infty,\infty}^{-s})}\big). 
\end{align}
Summing \eqref{0415-6}, \eqref{0415-7} and \eqref{0301-27}, we get 
\begin{align*}
  &{\|\nabla^\perp\Lambda^{-1}(\psi*w)\|}_{L^\infty(0,T;L^\infty)}
  +{\|w\|}_{L^\infty(0,T;B_{\infty,\infty}^{-s})}\notag\\
  \leq &C\Bigg({\|e^{-t\Lambda}\theta_0\|}_{\widetilde{L}^r(0,T;B_{\infty,\infty}^{1/r})}
  +\sum_{i=1}^{2}\Big(T{\|\theta^{(i)}\|}_{L^\infty(0,T;B_{\infty,\infty}^0)}
  +{\|\widetilde{\theta}^{(i)}\|}_{L^\infty(0,T;B_{\infty,\infty}^0)}\Big)\Bigg)\\
  &\times\big({\|\nabla^\perp\Lambda^{-1}(\psi*w)\|}_{L^\infty(0,T;L^\infty)}
  +{\|w\|}_{L^\infty(0,T;B_{\infty,\infty}^{-s})}\big). 
\end{align*}
Now from the assumption Theorem \ref{0308-1} \eqref{0308-4}, using Lemma \ref{0326-10} and Lemma \ref{0326-4}, 
we have 
\begin{equation*}
  \lim_{T\to 0}{\|e^{-t\Lambda}\theta_0\|}_{\widetilde{L}^r(0,T;B_{\infty,\infty}^{1/r})}=0
  \text{ and }
  \lim_{T\to 0}{\|\widetilde{\theta}^{(i)}\|}_{L^\infty(0,T;B_{\infty,\infty}^{0})}=0. 
\end{equation*}
Thus, due to a similar argument as for $3/2<\alpha\leq 2$, 
we can prove that the uniqueness of the mild solution of \eqref{SQG}. 

\subsection{Case $1<\alpha\leq 3/2$} 
Suppose $\theta^{(1)}, \theta^{(2)}\in C([0,T];B_{\infty,1}^{-1/2})$ (if $\alpha=3/2$) and 
$\theta^{(1)}, \theta^{(2)}\in C([0,T];B_{\infty,\infty}^{1-\alpha})$ (if $1<\alpha<3/2$). 
We estimate $w$ using $L^\infty(0,T;B_{\infty,\infty}^{1-\alpha})$ norm. 
By a similar argument to the case when $\alpha=1$, 
for $\alpha=3/2$, we obtain 
\begin{align*}
  &{\|\nabla^\perp\Lambda^{-1}(\psi*w)\|}_{L^\infty(0,T;L^\infty)}
  +{\|w\|}_{L^\infty(0,T;B_{\infty,\infty}^{-1/2})}\notag\\
  \leq &C\Bigg({\|e^{-t\Lambda^{\frac{3}{2}}}\theta_0\|}_{\widetilde{L}^r(0,T;B_{\infty,\infty}^{-1/2+3/2r})}
  +\sum_{i=1}^{2}\Big(T{\|\theta^{(i)}\|}_{L^\infty(0,T;B_{\infty,1}^{-1/2})}
  +{\|\widetilde{\theta}^{(i)}\|}_{L^\infty(0,T;B_{\infty,1}^{-1/2})}\Big)\Bigg)\\
  &\times\big({\|\nabla^\perp\Lambda^{-1}(\psi*w)\|}_{L^\infty(0,T;L^\infty)}
  +{\|w\|}_{L^\infty(0,T;B_{\infty,\infty}^{-1/2})}\big), 
\end{align*}
and for $1<\alpha<3/2$, we obtain 
\begin{align*}
  &{\|\nabla^\perp\Lambda^{-1}(\psi*w)\|}_{L^\infty(0,T;L^\infty)}
  +{\|w\|}_{L^\infty(0,T;B_{\infty,\infty}^{1-\alpha})}\notag\\
  \leq &C\Bigg({\|e^{-t\Lambda^\alpha}\theta_0\|}_{\widetilde{L}^r(0,T;B_{\infty,\infty}^{1-\alpha+\alpha/r})}
  +\sum_{i=1}^{2}\Big(T{\|\theta^{(i)}\|}_{L^\infty(0,T;B_{\infty,\infty}^{1-\alpha})}
  +{\|\widetilde{\theta}^{(i)}\|}_{L^\infty(0,T;B_{\infty,\infty}^{1-\alpha})}\Big)\Bigg)\\
  &\times\big({\|\nabla^\perp\Lambda^{-1}(\psi*w)\|}_{L^\infty(0,T;L^\infty)}
  +{\|w\|}_{L^\infty(0,T;B_{\infty,\infty}^{1-\alpha})}\big), 
\end{align*}
which implies the uniqueness of the mild solution of \eqref{SQG}. 

\subsection{Case $0<\alpha<1$} 
Now $\theta^{(1)},\theta^{(2)}\in C([0,T];B_{\infty,\infty}^{1-\alpha})$. 
We choose $s\in\mathbb{R}$ to satisfy $0<s<1-\alpha$. 
We estimate $w$ using $L^\infty(0,T;B_{\infty,\infty}^{-s})$ norm. 
In this case, we do not need to use the stracture of the nonlinear term of \eqref{SQG} (\eqref{0301-31} and Lemma \ref{0325-5}). 
Using the maximum principle (Lemma \ref{0326-7}), 
the product estimate (Lemma \ref{0325-8}), 
and the commutator estimate (Lemma \ref{0325-7} and Lemma \ref{0326-9}), we obtain 
\begin{align*}
  &{\|\nabla^\perp\Lambda^{-1}(\psi*w)\|}_{L^\infty(0,T;L^\infty)}
  +{\|w\|}_{L^\infty(0,T;B_{\infty,\infty}^{-s})}\notag\\
  \leq &C\Bigg({\|e^{-t\Lambda^\alpha}\theta_0\|}_{\widetilde{L}^r(0,T;B_{\infty,\infty}^{1-\alpha+\alpha/r})}
  +\sum_{i=1}^{2}\Big(T{\|\theta^{(i)}\|}_{L^\infty(0,T;B_{\infty,\infty}^{1-\alpha})}
  +{\|\widetilde{\theta}^{(i)}\|}_{L^\infty(0,T;B_{\infty,\infty}^{1-\alpha})}\Big)\Bigg)\\
  &\times\big({\|\nabla^\perp\Lambda^{-1}(\psi*w)\|}_{L^\infty(0,T;L^\infty)}
  +{\|w\|}_{L^\infty(0,T;B_{\infty,\infty}^{-s})}\big), 
\end{align*}
and it shows that ${\|w\|}_{L^\infty(0,T_0;B_{\infty,\infty}^{-s})}=0$ for some small $T_0>0$. 

\appendix 
\section{The product and the negative regularity }\label{0403-2}
In this section, we discuss the justification of the product $(\nabla^\perp\Lambda^{-1}f)g$ 
and the sum of two products $(\nabla^\perp\Lambda^{-1}f)g+(\nabla^\perp\Lambda^{-1}g)f$
for $f,g$ in the non-homogenous Besov spaces with the negative regularity indices. 
By Lemma \ref{0325-3}, if $s<0$, then for any $f,g\in B_{p,q}^s$, we have 
\begin{equation*}
  \sum_{|k-l|>1}(\phi_k*f)(\phi_l*g)\in \mathcal{S}', 
\end{equation*}
where the frequency of $(\phi_k*f)(\phi_l*g)$ is controled by $\phi_k*f$ or $\phi_l*g$ due to $|k-l|>1$. 
On the other hund, the following proposition holds. 
\begin{prop}\label{0404-3}
  Let $s\in\mathbb{R}$ and $1\leq p,q\leq\infty$. 
  If $s<0$, then there exist $f,g\in B_{p,q}^s(\mathbb{R}^2)$ such that 
  \begin{equation*}
    \sum_{|k-l|\leq 1}\left(\nabla^\perp\Lambda^{-1}(\phi_k*f)\right)(\phi_l*g)\notin\mathcal{S}' 
  \end{equation*}
  in the sense that there exists $\varphi\in\mathcal{S}$ such that 
  \begin{equation*}
    \sum_{|k-l|\leq 1}\int_{\mathbb{R}^2}\left(\partial_{x_j}\Lambda^{-1}(\phi_k*f)\right)(x)(\phi_l*g)(x)\varphi(x)~{\rm d}x=\infty,\ j=1,2.
  \end{equation*}
\end{prop}
\begin{proof}
  We prove only the case when $j=1$. 
  Let $\chi\in C_c^\infty(\mathbb{R}^2)$ with $\text{supp }\chi\subset\{\xi\in\mathbb{R}^2|\ |\xi|\leq 1/10\}$, 
  $\chi\geq 0$ and $\chi\equiv1$ near $0$. 
  We set $f$ and $g$ as follows: 
  \begin{equation*}
    f(x):=\lim_{N\to \infty}f_N(x)
    :=\lim_{N\to \infty}\sum_{n=1}^{N}2^{-sn}n^{-2}\mathcal{F}^{-1}\left[\chi(\xi-2^ne_1)\right](x)
  \end{equation*}
  and 
  \begin{equation*}
    g(x):=\lim_{N\to\infty}g_N(x)
    :=\lim_{N\to\infty}\sum_{n=1}^{N}2^{-sn}n^{-2}\mathcal{F}^{-1}\left[\chi(\xi+2^ne_1)\right](x), 
  \end{equation*}
  where $e_1=(1,0)$. 
  We have 
  \begin{align*}
    {\|f\|}_{B_{p,q}^s}
    \leq {\|f\|}_{B_{p,1}^s}
    \leq C\sum_{j\in\mathbb{N}}j^{-2}<\infty, 
  \end{align*}
  which implies $f\in B_{p,q}^s$. 
  Similariy, we get $g\in B_{p,q}^s$. 
  Next, for each $N\in\mathbb{N}$, we consider 
  \begin{equation*}
    \sum_{|k-l|\leq 1}\int_{\mathbb{R}^2}\left(\partial_{x_1}\Lambda^{-1}(\phi_k*f_N)\right)(x)(\phi_l*g_N)(x)\varphi(x)~{\rm d}x. 
  \end{equation*}
  We choose $\varphi=\widehat{\chi}(-\cdot)$. 
  By Plancherel's Theorem, we write 
  \begin{align*}
    &\int_{\mathbb{R}^2}\left(\partial_{x_1}\Lambda^{-1}(\phi_k*f_N)\right)(x)(\phi_l*g_N)(x)\widehat{\chi}(-x)~{\rm d}x\\ 
    =&\int_{\mathbb{R}^2}\int_{\mathbb{R}^2}\frac{i(\xi_1-\eta_1)}{|\xi-\eta|}\widehat{\phi}_k(\xi-\eta)\widehat{f}_N(\xi-\eta)\widehat{\phi}_l(\eta)\widehat{g}_N(\eta)~{\rm d}\eta\ \chi(\xi)~{\rm d}\xi, 
  \end{align*}
  and since $\chi$ has the support near the origin, we have 
  \begin{equation*}
    \chi(\xi)\widehat{f}_N(\xi-\eta)\widehat{g}_N(\eta)=\sum_{n=1}^{N}2^{-2sn}n^{-4}\chi(\xi)\chi(\xi-\eta-2^ne_1)\chi(\eta+2^ne_1). 
  \end{equation*}
  Also, since $\xi_1-\eta_1>0$ on $\text{supp}\chi(\cdot-2^ne_1)$ for every $n\in\mathbb{N}$, 
  we estimate from below by the sum with $k=l=n$. 
  Thus, we obtain 
  \begin{align*}
    &\sum_{|k-l|\leq 1}\int_{\mathbb{R}^2}\frac{\xi_1-\eta_1}{|\xi-\eta|}\widehat{\phi}_k(\xi-\eta)\widehat{f}_N(\xi-\eta)\widehat{\phi}_l(\eta)\widehat{g}_N(\eta)~{\rm d}\eta\ \chi(\xi)\\ 
    \geq &c\sum_{n=1}^{N}2^{-2sn}n^{-4}\int_{\mathbb{R}^2}\widehat{\phi}_n(\xi-\eta)\widehat{\phi}_n(\eta)\chi(\xi-\eta-2^ne_1)\chi(\eta+2^ne_1)~{\rm d}\eta\ \chi(\xi)\\ 
    \geq& c\sum_{n=1}^{N}2^{-2sn}n^{-4}, 
  \end{align*}
  where the constant $c>0$ is independent of $n$. 
  Since $s<0$, we have 
  \begin{equation*}
    \sum_{n=1}^{N}2^{-2sn}n^{-4}\to\infty, \text{ as } N\to\infty, 
  \end{equation*}
  which implies Proposition \ref{0404-3} holds. 
\end{proof}

On the other hund, we can define $(\nabla^\perp\Lambda^{-1}f)g+(\nabla^\perp\Lambda^{-1}g)f$ 
using the Fourier transform even in the case when $s<0$. 
Formally, if $f,g$ are smooth functions, then we have 
\begin{align*}
  &\sum_{|k-l|\leq 1}\left(\left(\nabla^\perp\Lambda^{-1}(\phi_k*f)\right)(\phi_l*g)
  +\left(\nabla^\perp\Lambda^{-1}(\phi_l*g)\right)(\phi_k*f)\right)\\ 
  =&\sum_{|k-l|\leq 1}\left(\left(\nabla^\perp\Lambda^{-1}(\phi_k*f)\right)(\phi_l*g)
  -\left(\Lambda^{-1}(\phi_l*g)\right)\left(\nabla^\perp (\phi_k*f)\right)\right)\\
  &+\sum_{|k-l|\leq 1}\nabla^\perp\left(\left(\Lambda^{-1}(\phi_l*g)\right)(\phi_k*f)\right)\\ 
  =&\sum_{|k-l|\leq 1}\nabla\cdot m(D_1,D_2)\left(\nabla^\perp\Lambda^{-1}(\phi_k*f),\Lambda^{-1}(\phi_l*g)\right)\\ 
  &+\sum_{|k-l|\leq 1}\nabla^\perp\left(\left(\Lambda^{-1}(\phi_l*g)\right)(\phi_k*f)\right), 
\end{align*}
where 
\begin{equation*}
  m(D_1,D_2)(f,g)
  :=\mathcal{F}^{-1}\left[\displaystyle\int_{\mathbb{R}^2}m(\eta,\xi-\eta)\widehat{f}(\eta)\widehat{g}(\xi-\eta)~{\rm d}\eta\right]
\end{equation*}
and 
\begin{equation*}
  m(\eta,\xi-\eta):=\displaystyle\frac{-i(\xi-2\eta)}{|\xi-\eta|+|\eta|}. 
\end{equation*}
Here, the bilinear Fourier multiplier $m(\cdot,\cdot)$ is bounded from $L^{p_1}\times L^{p_2}$ to $L^p$ ($1\leq p,p_1,p_2<\infty$ with $1/p=1/p_1+1/p_2$) (see \cite{Gr_Mi_To_2013}). 
Thus, if $f,g\in B_{p,q}^s$ with $s<0$ and $2\leq p<\infty$, 
then we define $\displaystyle\sum_{|k-l|\leq 1}\left(\left(\nabla^\perp\Lambda^{-1}(\phi_k*f)\right)(\phi_l*g)
+\left(\nabla^\perp\Lambda^{-1}(\phi_l*g)\right)(\phi_k*f)\right)$ as follows: 
\begin{dfn}\label{0415-8}
  Let $s<0$, $2\leq p<\infty$ and $1\leq q\leq \infty$. 
  For $f,g\in B_{p,q}^s(\mathbb{R}^2)$, we define 
  \begin{align*}
    &\subscripts{\mathcal{S}'}{\left\langle\sum_{|k-l|\leq 1}\left(\left(\nabla^\perp\Lambda^{-1}(\phi_k*f)\right)(\phi_l*g)
    +\left(\nabla^\perp\Lambda^{-1}(\phi_l*g)\right)(\phi_k*f)\right), \varphi \right\rangle}{\mathcal{S}}\\ 
    :=&\sum_{|k-l|\leq 1}\subscripts{\mathcal{S}'}{\left\langle\nabla\cdot m(D_1,D_2)\left(\nabla^\perp\Lambda^{-1}(\phi_k*f),\Lambda^{-1}(\phi_l*g)\right), \varphi\right\rangle}{\mathcal{S}}\\
    &+\sum_{|k-l|\leq 1}\subscripts{\mathcal{S}'}{\left\langle\nabla^\perp\left(\left(\Lambda^{-1}(\phi_l*g)\right)(\phi_k*f)\right), \varphi \right\rangle}{\mathcal{S}}. 
  \end{align*}
\end{dfn}

In the case of $p=\infty$, we consider the variables of $f$ and $g$ separately, 
then we can write 
\begin{align*}
  &\big(\nabla^\perp\Lambda^{-1}f\big)(x)g(x')
  -\big(\Lambda^{-1}g\big)(x')\big(\nabla^\perp f\big)(x') \notag\\
  =&\int_{\mathbb{R}^2}e^{ix\cdot\xi}\frac{i\xi^\perp}{|\xi|}\widehat{f}(\xi)~{\rm d}\xi
  \int_{\mathbb{R}^2}e^{ix'\cdot\eta}\widehat{g}(\eta)~{\rm d}\eta
  -\int_{\mathbb{R}^2}e^{ix'\cdot\eta}\frac{1}{|\eta|}\widehat{g}(\eta)~{\rm d}\eta
  \int_{\mathbb{R}^2}e^{ix\cdot\xi}i\xi^\perp\widehat{f}(\xi)~{\rm d}\xi\notag\\
  =&\nabla_{\mathbb{R}^4}\cdot \mathcal{F}^{-1}_{\mathbb{R}^4}\left[\frac{-i(-\xi,\eta)}{|\xi|+|\eta|}\frac{i\xi^\perp}{|\xi|}\widehat{f}(\xi)\frac{1}{|\eta|}\widehat{g}(\eta)\right](x,x'), 
\end{align*}
where $\mathcal{F}^{-1}_{\mathbb{R}^4}$ is the Fourier inverse transform in $\mathbb{R}^4$, 
and $\nabla_{\mathbb{R}^4}$ is the gradient in $\mathbb{R}^4$. 
Hence, if $p=\infty$, then we define $\displaystyle\sum_{|k-l|\leq 1}\left(\left(\nabla^\perp\Lambda^{-1}(\phi_k*f)\right)(\phi_l*g)
+\left(\nabla^\perp\Lambda^{-1}(\phi_l*g)\right)(\phi_k*f)\right)$ as well as Definition \ref{0415-8}. 

Also, if $\displaystyle\sum_{|k-l|\leq 1}\left(\left(\nabla^\perp\Lambda^{-1}(\phi_k*f)\right)(\phi_l*g)
+\left(\nabla^\perp\Lambda^{-1}(\phi_l*g)\right)(\phi_k*f)\right)\in\mathcal{S}'$, 
then we can define the sum of two products $(\nabla^\perp\Lambda^{-1}f)g+(\nabla^\perp\Lambda^{-1}g)f$. 
By Lemma \ref{0325-5}, 
we can justify the sum of two products $(\nabla^\perp\Lambda^{-1}f)g+(\nabla^\perp\Lambda^{-1}g)f$ when $s\geq -1/2$. 
\begin{lem}\label{0423-1}
  Let $-1/2<s<0$, $2\leq p\leq \infty$ and $1\leq q\leq\infty$ 
  or $s=-1/2$, $4\leq p\leq\infty$ and $1\leq q\leq p/(p-2)$. 
  If $p<\infty$, then for any $f,g\in B_{2p,q}^s$, we have 
  \begin{equation*}
    \sum_{|k-l|\leq 1}m(D_1,D_2)\left(\nabla^\perp\Lambda^{-1}(\phi_k*f),\Lambda^{-1}(\phi_l*g)\right)\in L^p(\mathbb{R}^2). 
  \end{equation*}
  Also, if $p=\infty$, then for any $f,g\in B_{\infty,q}^s$, we have 
  \begin{equation*}
    \sum_{|k-l|\leq 1}\mathcal{F}^{-1}_{\mathbb{R}^4}\left[\frac{-i(-\xi,\eta)}{|\xi|+|\eta|}\frac{i\xi^\perp}{|\xi|}\widehat{\phi}_k(\xi)\widehat{f}(\xi)\frac{1}{|\eta|}\widehat{\phi}_l(\eta)\widehat{g}(\eta)\right]\in L^\infty(\mathbb{R}^4). 
  \end{equation*}
\end{lem}

At the end of this section, 
we show an example where the product estimate in the Besov space does not hold 
for $(\Lambda^{-1}f)g$ when $s<-1/2$.

\begin{prop}\label{0405-1}
  Let $s,s'\in\mathbb{R}$, $2\leq p\leq\infty$ and $1\leq q\leq \infty$. 
  If $s<-1/2$, then there exist $f,g\in B_{2p,q}^s$ such that 
  \begin{equation*}
    {\|(\Lambda^{-1}f)g\|}_{B_{p,q}^{s'}}=\infty. 
  \end{equation*}
\end{prop}
\begin{proof}
  Let $\chi\in C_c^\infty(\mathbb{R}^2)$ with $\text{supp }\chi\subset\{\xi\in\mathbb{R}^2|\ |\xi|\leq 1/10\}$, 
  $\chi\equiv1$ near $0$ and $\chi\geq 0$. 
  For each $N\in\mathbb{N}$, we set 
  \begin{equation*}
    f_N(x)=g_N(x):=\sum_{n=1}^{N}2^{-sn}n^{-2}\mathcal{F}^{-1}\left[\chi(\xi-2^ne_1)+\chi(\xi+2^ne_1)\right](x)
  \end{equation*}
  and 
  \begin{equation*}
    f(x)=g(x):=\lim_{N\to\infty}f_N(x). 
  \end{equation*}
  By a similar argument to Proposition \eqref{0404-3}, we get $f,g\in B_{2p,q}^s$. 
  Also, using Young's inequality, we have 
  \begin{align*}
    {\|(\Lambda^{-1}f_N)g_N\|}_{B_{p,q}^{2s+1}}
    \geq& {\left\|\psi*\left((\Lambda^{-1}f_N)g_N\right)\right\|}_{L^{p}}\\ 
    \geq c&{\left\|\mathcal{F}^{-1}\left[\widehat{\psi}(\xi)\int_{\mathbb{R}^2}\frac{1}{|\xi-\eta|}\widehat{f}_N(\xi-\eta)\widehat{g}_N(\eta)~{\rm d}\eta\right]\right\|}_{L^\infty}\\ 
    \geq c&\int_{\mathbb{R}^2}\widehat{\psi}(\xi)\int_{\mathbb{R}^2}\frac{1}{|\xi-\eta|}\widehat{f}_N(\xi-\eta)\widehat{g}_N(\eta)~{\rm d}\eta~{\rm d}\xi. 
  \end{align*}
  Note that if $|\xi|\leq 4/3$, then for any $n,m\in\mathbb{N}$, we have 
  \begin{equation*}
    \chi(\xi-\eta-2^ne_1)\chi(\eta-2^me_1)
    =\chi(\xi-\eta+2^ne_1)\chi(\eta+2^me_1)
    =0. 
  \end{equation*}
  Thus, by a similar estimate to Proposition \eqref{0404-3}, we obtain 
  \begin{align*}
    &\int_{\mathbb{R}^2}\widehat{\psi}(\xi)\int_{\mathbb{R}^2}\frac{1}{|\xi-\eta|}\widehat{f}_N(\xi-\eta)\widehat{g}_N(\eta)~{\rm d}\eta~{\rm d}\xi\\ 
    \geq&\sum_{n=1}^N2^{(-2s-1)n}n^{-4} 
    \to \infty,\text{ as }N\to\infty, 
  \end{align*}
  which implies Proposition \ref{0405-1} holds. 
\end{proof}

\noindent
{\bf Data availability statement.} This manuscript has no associated data.\\
{\bf Conflict of Interest.} The author declares that he has no conflict of interest.

\end{document}